\documentclass[a4paper,10pt]{article}

\usepackage[round]{natbib}
\usepackage{xcolor}
\usepackage{amsmath}
\usepackage{amsthm}
\usepackage{amsfonts}
\usepackage{mathrsfs}
\usepackage{enumerate}
\usepackage{eurosym}

\author{ Vaios Laschos$^*$ \and Klaus Obermayer\thanks{Fakult\"at Elektrotechnik und Informatik, and Bernstein Center for Computational Neuroscience, Technische Universit\"at Berlin,  Marchstr.~23, 10587, Berlin, Germany ({\tt yun.shen3@gmail.com, vaios.laschos@tu-berlin.de, klaus.obermayer@mailbox.tu-berlin.de}).} \and Yun Shen$^*$  \and Wilhelm Stannat\thanks{Institut f\"ur Mathematik, and Bernstein Center for Computational Neuroscience, Technische Universit\"at Berlin, Stra{\ss}e des 17.\ Juni 136, 10623, Berlin, Germany ({\tt stannat@math.tu-berlin.de}).} }

\date{March 16, 2018}

\usepackage[margin=1.3in]{geometry}

\newtheorem{theorem}{\sc Theorem}[section]
\newtheorem{lemma}[theorem]{\sc Lemma}
\newtheorem{proposition}[theorem]{\sc Proposition}
\newtheorem{definition}[theorem]{\sc Definition}
\newtheorem{corollary}[theorem]{\sc Corollary}

\theoremstyle{remark}
\newtheorem{remark}[theorem]{{\it Remark}}
\newtheorem{assumption}[theorem]{{\it Assumption}}
\newtheorem{example}[theorem]{{\it Example}}
\newtheorem{algorithm}[theorem]{{\it Algorithm}}

\DeclareMathOperator*{\lip}{Lip}
\DeclareMathOperator*{\aes}{\textrm{\AE}}
\DeclareMathOperator*{\epi}{\textrm{epi}}
\DeclareMathOperator*{\dom}{\textrm{dom}}

\usepackage{enumerate}

\newcommand{\proofbegin}{\noindent {\sc Proof.}\ }
\newcommand{\proofend}{{\quad $\square$}}

\title{A Fenchel-Moreau-Rockafellar type theorem on the Kantorovich-Wasserstein space with Applications in Partially Observable Markov Decision Processes}

\begin{document}
 
\maketitle

\abstract{By using the fact that the space  of all probability measures with finite support can be somehow completed in two different fashions, one generating the Arens-Eells space and another generating the Kantorovich-Wasserstein (Wasserstein-1) space, and by exploiting the duality relationship between the Arens-Eells space with the space of Lipschitz functions, we provide a dual representation of Fenchel-Moreau-Rockafellar type for proper  convex functionals on Wasserstein-1. We retrieve dual transportation inequalities as a Corollary and we provide examples where the theorem can be used to easily prove dual expressions like the celebrated Donsker-Varadhan variational formula. Finally our result allows to write convex functions as the supremum over all linear functions that are generated by roots of its conjugate dual, something that we apply to the field of Partially observable Markov decision processes (POMDPs) to approximate the value function of a given POMDP by iterating level sets. This extends the method used in \citet{smallwood1973optimal} for finite state spaces to the case were the state space is a Polish metric space.}\newline



\section{Introduction}

Dual representation, of the Fenchel-Moreau-Rockafellar type (see e.g., \citet{bot2010,ioan2009duality,zalinescu2002convex}), plays an important role in convex analysis and has wide applications in various fields \citep{borwein2010convex,boyd2004convex}. In this paper, we are interested in a conjugate dual representation of functions on probability measures, for which the ``choice of dual space'' will allow for real-world applications in the field of Partially observable Markov decision processes (POMDPs) along the lines of \citet{smallwood1973optimal}. In that direction, we are going to use the special connection of the Wasserstein-1 space with the space of Lipschitz functions. 

More specifically, for $(\mathsf X, d)$ being a Polish space, $(\mathscr P_{1}(\mathsf X),W_{1})$ the Wasserstein-1 space over $\mathsf X,$ and $\mathscr L(\mathsf X)$ the space of all real-valued Lipschitz functions on $(\mathsf X, d)$, the \emph{conjugate function} $\rho: \mathscr L(\mathsf X) \rightarrow \bar{\mathbb R}$ of a function $\phi: \mathscr P_{1}(\mathsf X) \rightarrow \bar{\mathbb R}$ is defined as 

\begin{align}\label{eq:conjugate}
 \rho(f) := \sup_{\mu \in \mathscr P_{1}(\mathsf X)} \left( \int f d\mu - \phi(\mu) \right)
\end{align}
and the \emph{second conjugate function} $\phi^c: \mathscr P_{1}(\mathsf X) \rightarrow \bar{\mathbb R}$ is defined as
\begin{align}
 \phi^c(\mu) := \sup_{f \in \mathscr L(\mathsf X)} \left( \int f d \mu - \rho(f) \right).
\end{align}
A justification that $\rho$ is well defined and proper can be found in the proof of \cite[Theorem 2.3.5]{ioan2009duality}. The first main result of this paper is the following conjugate duality.

\begin{theorem}\label{th:fm}
 Let $\phi: \mathscr P_{1}(\mathsf X) \rightarrow \bar{\mathbb R}$ be a proper convex function on $(\mathscr P_{1}(\mathsf X), W_1),$ i.e. a convex and lower semicontinuous function, satisfying $\phi(\mu) > - \infty$ for all $\mu \in \mathscr P_{1}(\mathsf X)$ and $\phi(\mu_{0}) \in \mathbb R$ for some $\mu_{0} \in \mathscr P_{1}(\mathsf X)$. Then $\phi(\mu) = \phi^c(\mu)$, $\forall \mu \in \mathscr P_{1}(\mathsf X)$.
\end{theorem}

In \citet[Theorem 5.26]{villani2009optimal}, one has to assume that $\phi(\mu) = \phi^c(\mu),$  were the conjugates there are defined by taking the supremum over  $C_{b}(\mathsf X),$ in order to establish a pair of dual inequalities connecting an optimal transport distance to $\phi(\mu).$ Theorem \ref{th:fm} implies that, in the case of the Wasserstein-1 distance, this assumption is always satisfied, provided that $\mathscr L (\mathsf X)$ is used in place of  $C_{b}(\mathsf X)$ for defining the conjugate dual. We will continue by  providing a Wasserstein-1 specific version of \citet[Theorem 5.26]{villani2009optimal}, and an example of a pair of well known dual functions where is it very simple to calculate $\rho(\mu),$ but the known proofs of  $\phi^c=\phi,$ are technically more complicated.

\begin{corollary}\label{th:vil}
	Let $\phi: \mathscr P_{1}(\mathsf X) \rightarrow \bar{\mathbb R}$ be a proper convex function on $(\mathscr P_{1}(\mathsf X), W_1).$ Let also  $\rho$ be its conjugate as in \eqref{eq:conjugate}. Let finally $\Phi$ a real increasing and convex function with $\Phi(0)=0.$ We have
	\begin{equation*}
	\Phi( W_1(\mu, \nu) )\leq \phi(\mu),\, \forall \mu\in\mathcal{P}_{1}(\mathsf{X}) \Leftrightarrow \rho\left(\int_{X}tfd\nu-tf-\Phi^{*}(t)\right)\leq 0, \forall f\in [\mathscr L(\mathsf X)]_1, t\in\mathbb{R}
	\end{equation*} 
where $\Phi^{*}$ is the Legendre dual, i.e. is given by the formula $\Phi^{*}(s)=\sup\{st-\Phi(t)\},$ and $[\mathscr L(\mathsf X)]_1$ is the set of all Lipschitz functions with constant 1 (see next section for definition). 

\end{corollary}

The proof of the Corollary is straightforward and can be found in the Appendix.  We would like to remark that  $[\mathscr L(\mathsf X)]_1$ can be substituted with any subset $\mathcal{A}$ of $[\mathscr L(\mathsf X)]_1$ that satisfies $\mathcal{A}+\mathbb{R}=[\mathscr L(\mathsf X)]_1.$ We proceed now with an example.

\begin{example}
We will show that a direct application of Theorem \ref{th:fm} can provide  the celebrated Donsker-Varadhan variational formula, which has many fundamental applications in the theory of large deviations (see for example \cite{Dupuis1997}, where a whole class of large deviation principles are proved by applying the formula) and in statistical physics in general.

It is known that the following pair of dual equation hold:
\begin{equation}\label{first var}
\log\int_{\mathsf X}e^{g}d\nu=\sup_{\mu\in\mathcal{P}(\mathsf X)}\left\{ \int_{\mathsf X}gd\mu-\mathcal{R}(\mu|\nu)\right\}, \hspace{16pt}\forall g\in C_{b}(\mathsf X)
\end{equation}
and
\begin{equation}\label{secvar}\begin{split}
\mathcal{R}(\mu|\nu) = \sup_{g\in C_{b}(\mathsf X)}\left\{\int_{\mathsf X}gd\mu-\log\int_{\mathsf X}e^{g}d\nu\right\}, \hspace{16pt}\forall \mu\in\mathcal{P}(\mathsf{X})
\end{split}
\end{equation}
where $\mathcal{R}(\mu|\nu)$ is the relative entropy functional given by
\begin{equation*}
\mathcal{R}(\mu|\nu)=\begin{cases}\int_{\mathsf X}\frac{d\mu}{d\nu}\log\left(\frac{d\mu}{d\nu}\right)d\nu&\text{if}\hspace{8pt} \mu<<\nu,\\\infty&\text{otherwise}.\end{cases}
\end{equation*}
As it is shown in Lemma \ref{firstvar} in the Appendix, the first formula is straightforward to prove, even when $C_{b}(\mathsf X)$ is replaced by $\mathscr{L}(\mathsf X)$ (the original proof is even simpler and can be found in page 34 of \cite{Dupuis1997} ). On the other hand, the second requires a more technical proof as one can see in \cite[Lemma 1.4.3.]{Dupuis1997}. By applying Theorem \ref{th:fm}, we get the following alternative variational form for the relative entropy, i.e.:
\begin{equation}\label{secondvar2}
\hspace{-9pt}\mathcal{R}(\mu|\nu)=\sup_{g\in \mathscr L(\mathsf X)}\left\{\int_{\mathsf X}gd\mu-\log\int_{\mathsf X}e^{g}d\nu\right\}.
\end{equation}
If one wishes to further retrieve \eqref{secvar} from \eqref{secondvar2}, then it is a matter of a simple approximation argument as it is illustrated in Lemma \ref{con=lip}. Finally, applying Corollary \ref{th:vil} with $\Phi(t)=\frac{1}{2}t^{2}$ on our example, we can retrieve Bobkov and G\"{o}tze theorem \cite[Theorem 1.3]{bobkov}, i.e.
\begin{equation*}
\begin{split}
\int_{X}e^{tf}d\nu\leq e^{ct^{2}/2},&\, \hspace{8pt}\forall f\in \left\{f\in[\mathscr L(\mathsf X)]_1 \Big| \int_{\mathsf{X}}fd\nu=0\right\},\,t\in\mathbb{R} \hspace{6pt} \Leftrightarrow \\ &W_1(\mu, \nu)\leq \frac{1}{c}\sqrt{\mathcal{R}(\mu|\nu)},\hspace{8pt}\forall\mu\in\mathcal{P}_{1}(\mathsf{X}).
\end{split}
\end{equation*}
One could probably derive more refined versions of the inequality by playing around with the choice of $\Phi.$
\end{example}

Another fundamental application of Theorem \ref{th:fm} is to derive a conjugate dual form for the optimality equation for POMDPs. In order to do that we first are going to show that each function that satisfies the condition in Theorem \ref{th:fm} has a representation as a supremum over a suitable class of linear functionals. More specifically we have the following.

Let $\phi$ be a function on $\mathscr P_{1}(\mathsf X)$ and $\rho$ be its conjugate as in \eqref{eq:conjugate}. Consider the following sets 
\begin{align*}
\mathcal N_\phi := & \left\{ f \in \mathscr L(\mathsf X) \middle| \rho(f) = 0 \right\}, \textrm{ and } \\
\bar{\mathcal N}_\phi := & \left\{ f \in \mathscr L(\mathsf X) \middle| \rho(f) \leq 0 \right\} = \left\{ f \in \mathscr L(\mathsf X) \middle| \int f d \mu \leq \phi(\mu), \forall \mu \in \mathscr P_{1}(\mathsf X) \right\} .
\end{align*}
We call $\mathcal N_\phi$ the \emph{null level-set of $\phi$}, whereas the latter set $\bar{\mathcal N}_\phi$ is called the \emph{acceptance set of $\phi$} (cf.\ \cite[Section 4.1]{follmer2004stochastic}). Note that since $\rho$ is convex and lower semicontinuous, $\bar{\mathcal N}_\phi$ is convex and closed (see e.g.\ \cite[Theorem 2.2.9]{ioan2009duality}). %To see this, extending the domain of $\phi$ to contain the 0 measure such that $\phi(0) := 0$ yields $\phi(\mu) - \phi(0) \geq \int f (d\mu - d0), \forall \mu$, for any $f \in \bar{\mathcal N}_\phi.$ 
We have the following dual representation.
\begin{corollary}\label{coro:fm}
	Let $\phi: \mathscr P_{1}(\mathsf X) \rightarrow \bar{\mathbb R}$ be a function satisfying the condition in Theorem \ref{th:fm}. Then $$\phi(\mu) = \sup_{f \in \mathcal N_\phi} \int f d\mu =  \sup_{f \in \bar{\mathcal N}_\phi} \int f d\mu.$$
\end{corollary}

The proof is given in the end of the next section. Before we proceed we will provide some background on POMDPS and we are going to explain how  Corollary \ref{coro:fm} is meant to be used in that setting.

 A POMDP is a tupe of contrroled stochastic processes where it is assumed that the system dynamics are determined by a \emph{Markov process}, but the agent cannot directly observe the underlying state. POMDPs have important applications in various fields, such as operations research \citep{lovejoy1991survey}, robotics \citep{pineau2006anytime} and artificial intelligence \citep{kaelbling1998planning}. It is known that a POMDP can be reduced to a standard \emph{Markov decision process} (MDP) by using appropriate probability distributions over the hidden states. We refer to \citet{sondik1978optimal} for finite spaces, to \citet{sawaragi1970discrete} for countable spaces, and to \citet[Chapter 4]{hernandez1989adaptive} and \citet{feinberg2016partially} for Borel spaces. 

In the setting of finite state spaces with a discounted infinite-horizon objective, the value function encoding the maximum reward can be found by solving the following equation
\begin{align}
 \phi(\mu) = \max_{a} \!\left\{ \tilde{r}(\mu, a) + \alpha \sum_{y} \phi(\mu'(\mu, a, y)) \left( \sum_{x'}  P(x'|\mu, a) Q(y|x', a) \right) \right\}, \forall \mu \in \mathscr P(\mathsf X). \label{eq:phi}
\end{align}
Here, $x'$, $a$ and $y$ denote the hidden (or latent) state, the action and the observation, respectively, and $\mu \in \mathscr P(\mathsf X)$ is the distribution over states, while $\mu'$ is the posterior distribution of the successive state given by
\begin{align*}
\mu'(\cdot|\mu, a, y) = \frac{P(\cdot|\mu, a) Q(y|\cdot, a)}{\sum_{x'}  P(x'|\mu, a) Q(y|x', a)}.
\end{align*}
$P$ controls the transition probability between states, $Q$ models the observation probability of $y$ given states and actions, $r$ denotes the reward function, and $\alpha \in (0,1)$ serves as a discount factor. A formal introduction of POMDPs on Borel spaces can be found in Section \ref{sec:pomdps}. 

In spite of knowing the existence of such a theoretical solution, POMDPs were notoriously difficult to solve in practice \citep{shani2013survey}. In the case where the underlying space is finite, a fast algorithm called \emph{point-based value iteration} \citep{pineau2006anytime} was designed to overcome this numerical difficulty. This algorithm is mainly based on the property first observed by \citet{smallwood1973optimal}, that the optimal solution to equation \eqref{eq:phi} can be arbitrarily well approximated by a function of the following \emph{dual representation}
\begin{align}
 \phi(\mu) = \max_{h \in N} \sum_{x} h(x) \mu(x) \label{eq:dual},
 \end{align}
where $N$ can be chosen to be a \textbf{finite} collection of real functions on the hidden state space.

To our best knowledge, it is still an open question whether a similar approximation is also possible for POMDPs on continuous state spaces. The second major contribution of this paper is to provide an affirmative answer to this open question (see Theorem \ref{Secondbigtheorem}), by applying the conjugate duality on Wasserstein-1 spaces that we obtained in Theorem \ref{th:fm}. We would like to remark, that contrary to the finite state space case, the set $N$ appearing in the counterpart of \eqref{eq:dual} is uncountable, so a computable algorithm is still elusive, and it remains an \textbf{open problem} to find one. However, new advances in the field of neural networks involving the Wasserstein-1 distance, made substitution of collections of Lipschitz functions by appropriate finite sets a necessity \citep{Arjovsky}, and we are planning to perform further research in the immediate future.

\paragraph{Preliminaries} A \emph{Polish space}  $(\mathsf X,d)$ is a complete separable metric space and a \emph{Borel space} is a Borel subset of a Polish space. We denote by $C_{b}(\mathsf X)$ the set of all real bounded continuous functions on $\mathsf X.$   We further denote  the set of all probability measures on $\mathsf X$ by $\mathscr{P}(\mathsf X)$. Furthermore, we define the set
\begin{equation}\label{p1}
\mathscr P_{1}(\mathsf X) := \left\{ \mu \in  \mathscr P(\mathsf X) \ \middle| \ \int d(x_0, x) \mu(dx) < \infty \right\},
\end{equation}
where $x_0 \in \mathsf X$ is arbitrary. If $\mathsf Y$ is a Borel space, its Borel $\sigma$-algebra is denoted by $\mathcal B(\mathsf Y)$. 

%Let $\mathsf Y$ and $\mathsf Z$ be two Borel spaces. A \emph{stochastic kernel on $\mathsf Y$ given $\mathsf Z$} is a function $\psi(B|z), B \in \mathcal B(\mathsf Y), z \in \mathsf Z$ such that (i) $\psi(\cdot|z)$ is a probability measure on $\mathcal B(\mathsf Y)$ for every fixed $z \in \mathsf Z$, and (ii) $\psi(B|\cdot)$ is a measurable function on $\mathsf Z$ for every fixed $B \in \mathcal B(\mathsf Y)$. %For a Borel space $\mathsf X$, denote by $\mathscr P(\mathsf X)$ the space of all probability measures on $(\mathsf X, \mathcal B(\mathsf X))$.

\section{A conjugate duality on the Wasserstein-1 space}

\subsection{The Wasserstein-1 space}
 Let $(\mathsf X, d)$ be a Polish space with metric $d$. The Wasserstein-1 space over  $(\mathsf X, d)$ is the set $\mathscr{P}_{1}(\mathsf X),$ defined in \eqref{p1}, equipped with the \emph{Wasserstein-1 metric} given by
\begin{align}
 W_1(\mu, \nu) := \inf \left\{ \mathbb E \left[ d(X, Y) \right] \ \middle| \ \textrm{law}(X) = \mu, \textrm{law}(Y) = \nu \right\}. \label{eq:wmetric}
\end{align}
For any real-valued function $f: \mathsf X \rightarrow \mathbb R$, its Lipschitz seminorm is defined as 
\begin{align*}
 \lVert f \rVert_{\lip} := \sup_{x \neq y} \frac{\lvert f(x) - f(y) \rvert}{d(x,y)}.
\end{align*}
Denote by $\mathscr L(\mathsf X):= \left\{ f: \mathsf X \rightarrow \mathbb R \hspace{2pt}\middle|\hspace{2pt} \lVert f \rVert_{\lip} < \infty \right\}$ the space of all real-valued Lipschitz functions on $\mathsf X$ and by $[\mathscr L(\mathsf X)]_1 := \{ f\in \mathscr L \hspace{2pt}|\hspace{2pt} \lVert f \rVert_{\lip}\leq 1 \}$ the unit ball of $\mathscr L(\mathsf X)$. Then, it can be shown  \citep[Chapter 5]{villani2009optimal} that $W_1$ has the following representation:
\begin{align}
 W_1(\mu,  \nu) = \sup_{f \in [\mathscr L(\mathsf X)]_1} \left( \int f d\mu - \int f d\nu \right). \label{eq:kr-distance}
\end{align}

  We have the following property \citep[Theorem 6.18]{villani2009optimal}: {\it If $(\mathsf X, d)$ is Polish, $(\mathscr P_{1}(\mathsf X), W_1)$ is also Polish}.
  In the sequel, we are going to make use of the following set  $\mathscr D(\mathsf X)$ of all probability measures with finite support. More specifically, we define
\begin{align}
 \mathscr D(\mathsf X) := \left\{ \sum_{i=1}^n a_i \delta_{x_i} \middle| n \in \mathbb N, a_i \in \mathbb R_+, x_i\in\mathsf X, i=1,2, \ldots, n, \sum_{i=1}^n a_i = 1 \right\}. \label{eq:D}
\end{align}
By Theorem 6.18 in \cite{villani2009optimal}, $\mathscr D(\mathsf X)$ is dense in $(\mathscr P_{1}(\mathsf X), W_1)$.

\subsection{The Arens-Eells space} We recall some results of the Arens-Eells space based on \citet[Section 2.2 and 2.3]{weaver1999lipschitz}. 

\begin{definition}\it
 Let $(\mathsf X, d)$ be a metric space. A \emph{molecule} of $\mathsf X$ is a function $m: \mathsf X \rightarrow \mathbb R$ which is supported on a finite set and which satisfies $\sum_{x \in \mathsf X} m(x) = 0.$
\end{definition}

For any $x,y \in \mathsf X$ define the molecule $m_{xy} = \mathbf 1_x - \mathbf 1_y$, where $\mathbf 1_x$ denotes the indicator function on the singleton set $\{x\}$. Define the following seminorm for every molecule $m$:
\begin{equation}\label{aes}
 \lVert m \rVert_{\aes(\mathsf X)} := \inf \left\{ \sum_{i=1}^n \lvert a_i \rvert d(x_i, y_i) \middle| m= \sum_{i=1}^n a_i m_{x_i y_i} \right\},
\end{equation}
and let $\aes(\mathsf X)$ be the completion of the space of molecules, which is also called \emph{Arens-Eells space}. 

The following theorem states that $\aes(\mathsf X)$ is a predual of the space of Lipschitz functions  $\mathscr L(\mathsf X).$ This predual is unique in many important cases \citep{Weaver2016}. 
\begin{theorem}\citep[Theorem 2.2.2]{weaver1999lipschitz}\label{thm:aes}
	Let $(\mathsf X, d)$ be a metric space with at least one point $x_0$. Then  $(\aes(\mathsf X),\|\cdot\|_{\aes})^{*} \cong (\mathscr L(\mathsf X),\|\cdot\|_{\lip}).$
\end{theorem}

For any $f \in \mathscr L(\mathsf X)$ and $m \in \aes(\mathsf X)$, we define
$$\langle f, m \rangle: = \sum_{x \in \mathsf X} f(x) m(x).$$

\begin{corollary}\label{coro:aes}
	(i) $\lVert m \rVert_{\aes(\mathsf X)} = \max_{f \in [\mathscr L(\mathsf X)]_1} \langle f, m \rangle$. (ii) $\lVert \cdot \rVert_{\aes(\mathsf X)}$ is a norm on $\aes.$ 
\end{corollary}
For the proof, see \citet[Corollary 2.2.3]{weaver1999lipschitz}. \newline

 Replacing the Dirac measure in $\mathscr D(\mathsf X) $ by the indicator function, we define the following set of real-valued functions on $\mathsf X$
\begin{align}
\mathcal D(\mathsf X) := \left\{ \sum_{i=1}^n a_i \mathbf 1_{x_i} \middle| n \in \mathbb N, a_i \in \mathbb R_+, i=1,2, x\in \mathsf X, \ldots, n, \sum_{i=1}^n a_i = 1 \right\}. \label{eq:Dtilde}
\end{align}
Let $\psi: \mathscr D(\mathsf X) \rightarrow \mathcal D(\mathsf X)$ be defined by
\begin{equation}
\psi(\nu)(x):=\nu(\{x\})
\end{equation}

Obviously, $\psi$ is a bijection between $\mathscr D(\mathsf X)$ and $ \mathcal D(\mathsf X).$ Let now $\Psi: \mathscr D(\mathsf X)\times \mathscr D(\mathsf X)\rightarrow \aes(\mathsf X)$ be defined by  
\begin{align}
\Psi(\nu,\nu_{0}):= \psi(\nu) - \psi(\nu_0). 
\end{align}

To see that $\Psi(\nu,\nu_{0})$ is actually an element of $\aes(\mathsf X),$  notice that 

\begin{align*}
\sum_{x \in \mathsf X} \Psi(\nu,\nu_0)(x)= \sum_{x \in \mathsf X} \left(\psi(\nu)(x) - \psi(\nu_0)(x) \right)= 0, \forall \nu,\nu_{0} \in \mathscr D(\mathsf X).
\end{align*}

We remark that for every $\nu_{0}\in\mathscr D(\mathsf X),\hspace{2pt} \Psi$ is an injection from $\mathscr D(\mathsf X)\times\{\nu_{0}\}$ into  $\aes(\mathsf X).$ However, this obviously is not true for the whole product $\mathscr D(\mathsf X)\times \mathscr D(\mathsf X).$ Furthermore, $\Psi$ is a surjection into the set of molecules with ``total mass'' $\sum_{x\in\mathsf X}|m(x)|$ equal or less than $2$.

\subsection{Connecting Arens-Eells and Wasserstein-1 spaces.}
Recall that $\mathscr D(\mathsf X)$ is the subset of  $\mathscr P(\mathsf X),$ that contains all probability measures with finite support, and $\mathcal D(\mathsf X)$ defined in \eqref{eq:Dtilde} is its corresponding space of functions on $\mathsf X$, with $\psi: \mathscr D(\mathsf X) \rightarrow \mathcal D(\mathsf X)$ being the bijective map. 

\begin{proposition}\label{asas}
	$W_1(\nu, \nu_0) =  \lVert \Psi(\nu,\nu_0) \rVert_{\aes(\mathsf X)}=\lVert \psi(\nu) - \psi(\nu_0) \rVert_{\aes(\mathsf X)}$, for every $\nu,\nu_{0} \in \mathscr D(\mathsf X).$ 
\end{proposition}
\proofbegin 
By the dual representation \eqref{eq:kr-distance} of $W_1$, we have
\begin{align*}
W_1(\nu, \nu_0) = & \sup_{f \in [\mathscr L(\mathsf X)]_1} \left(\int f d \nu - \int f \nu_0\right) = \sup_{f \in [\mathscr L(\mathsf X)]_1} \sum_{x \in \mathsf X} \left( \psi(\nu)(x) - \psi(\nu_0)(x) \right)f(x)\\& = \lVert \psi(\nu) - \psi(\nu_0) \rVert_{\aes(\mathsf X)}=\lVert \Psi(\nu,\nu_0) \rVert_{\aes(\mathsf X)},
\end{align*}
where the second to the last equality is due to Corollary \ref{coro:aes}(i).
\proofend

Before we proceed, we would like to highlight the connection of the Wasserstein-1 space with the Arens-Eells space and its dual, namely the space of Lipschitz functions. We saw in the previous subsection and in Proposition \ref{asas} that we can embed the set $\mathscr D(\mathsf X)\times\mathscr D(\mathsf X)$ of pairs $(\nu,\nu_{0}) \in \mathscr P_{1}(\mathsf X)\times \mathscr P_{1}(\mathsf X)$ in the vector space $(\aes(\mathsf X), \lVert \cdot \rVert_{\aes(\mathsf X)})$, in a way that the Wasserstein-1 distance  $ W_1(\nu, \nu_0)$ is equal to the norm of the vector $\lVert \Psi(\nu,\nu_0) \rVert_{\aes(\mathsf X)}.$ Moreover, for every element $m$ of the set of molecules (which is dense in $\aes(\mathsf X)$), one can find a pair $(\nu^{(m)},\nu_{0}^{(m)})\in \mathscr D(\mathsf X)\times\mathscr D(\mathsf X),$ and a positive number $a^{(m)},$ such that $m=a^{(m)}\Psi(\nu^{(m)},\nu^{(m)}_{0}),$ and $\|m\|_{\aes(\mathsf X)}=a^{(m)} W_1(\nu^{(m)}, \nu^{(m)}_0).$ One can picture $\mathscr D(\mathsf X)\times\mathscr D(\mathsf X)$ as an absorbing set \citep{Schaefer1971} of a dense subspace of $\aes(\mathsf X),$ which hints that  $\mathscr D(\mathsf X)\times\mathscr D(\mathsf X)$ may be enough for characterizing its dual space $\mathscr L(\mathsf X)$.
 
  In the case where $\mathsf X$ is compact, one can get a better intuition by reading the exposition on the relation between the so-called Kantorovich-Rubenstein space, the Arens-Eells space and the space of Lipschitz functions in \citet[Chapter 2, Section 3]{weaver1999lipschitz} or \citet[Section VIII.4]{KantAki1982}.
  
   In addition, due to \eqref{eq:kr-distance}, the Wasserstein-1 space can be considered as a subspace of the dual of  $\mathscr L(\mathsf X).$ In what follows, we are going to exploit these relationships to prove a separation theorem on the Wasserstein space, and then our first main  result, Theorem \ref{th:fm}.

\textbf{Open Problem}.  It is tempting to generalize our approach to establish duality w.r.t. the Wasserstein-p distance with $p>1$.  To this end, one could try to generalize the Arens-Eells space to some ``p-version", for example by using $d(x_{i},y_{i})^{p}$ instead of $d(x_{i},y_{i})$ in \eqref{aes}. However, in this case, the  seminorm constructed in an analogous way to \eqref{aes}, is equal to zero everywhere. The approach we are following in order to prove the duality result, heavily depends on the fact that adding the same finite measure on two measures $\nu,\nu_{0}$ does not change its distance; here, one can understand the Wasserstein distance not between probability measures anymore, but between measures with the same total mass. Therefore, it remains an open problem, how one can generalize the result for Wasserstein-$p$ spaces with $p>1$. 
 
\subsection{A separation theorem on $\mathscr D(\mathsf X)$}  We are now ready to state a separation theorem on $\mathscr D(\mathsf X)$. We first restate the Nirenberg-Luenberger theorem (see \citet[Section 5.13]{luenbergeroptimization}, also known as the \emph{minimum norm duality theorem}).

\begin{definition}\it
 Let $K$ be a convex set in a real normed vector space $X$. Let $X^*$ be its dual space. The function $h_{K}(x^*) = \sup_{x \in K} \langle x, x^* \rangle$ on $X^*$ is called the support functional of $K$.
\end{definition}

\begin{theorem}\label{thm:nlt}
 Let $x_1$ be a point in a real normed space $X$ and let $d > 0$ denote its distance from the convex set $K$ having support functional $h_{K}$. Then
 \begin{align*}
  d = \inf_{x \in K} \lVert x - x_1 \rVert = \max_{\lVert x^* \rVert \leq 1} \left[ \langle x_1, x^*\rangle - h_{K}(x^*) \right]
 \end{align*}
 where the maximum on the right is attained by some $x_0^* \in X^*$.
\end{theorem}

\begin{theorem}\label{thm:mnd}
 Let $K$ be a convex set of $\mathscr D(\mathsf X)$ and $\nu_0 \in \mathscr D(\mathsf X)$ be a point not contained in $K$ satisfying $W_1(\nu, \nu_0) \geq \epsilon > 0, \forall \nu \in K$. Then, there exists a function $f \in [\mathscr L(\mathsf X)]_1$ such that $\int f d\nu \geq \int f d\nu_0 + \epsilon, \forall \nu \in K.$
\end{theorem}
\proofbegin
 Since $K$ is convex, $K' := \{\psi(\nu) - \psi(\nu_0) \mid \nu \in K \}$ is also convex in $\aes(\mathsf X).$ 
 To apply Theorem \ref{thm:nlt}, we set $X = \aes(\mathsf X)$ and its dual is $X^* = \aes^*(\mathsf X) = \mathscr L(\mathsf X)$. By the definition of $h_{K'}$, we have
 \begin{align*}
  h_{K'}(f) = \sup_{m \in K'} \langle f, m \rangle = \sup_{\nu \in K} \left\{ \int f d\nu - \int f  d\nu_0 \right\}, \forall f \in \mathscr L(\mathsf X).
 \end{align*}
 Then, by Theorem \ref{thm:nlt} and Proposition \ref{asas}, we obtain
 $$
  \epsilon \leq \inf_{\nu \in K} W_1(\nu, \nu_0) = \inf_{m \in K'} \lVert m \rVert_{\aes} = \max_{f \in [\mathscr L(\mathsf X)]_1 } \left[ - h_{K'}(f) \right],
 $$
which yields
$
 \max_{f \in [\mathscr L(\mathsf X)]_1} \left[ \inf_{\nu \in K} \int (-f) d \nu - \int (-f) d\nu_0 \right] \geq \epsilon.
$
Suppose the maximum is attained at $f_0$. Then
$
 \inf_{\nu \in K} \int (-f_0) d \nu - \int (-f_0) d\nu_0 \geq \epsilon,
$
which implies $\int (-f_0) d \nu \geq \int (-f_0) d\nu_0 + \epsilon, \forall \nu \in K.$
\proofend

In the following sections, we consider the dual space of the Cartesian product $\aes(\mathsf X) \times \mathbb R$, where we can obtain a similar result as above. The operator $\vee$ is defined as $a \vee b := \max(a,b).$
\begin{theorem}
 Let $\tilde K$ be a convex subset of $\mathscr D(\mathsf X) \times \mathbb R$ and $(\nu_0, r_0) \in \mathscr D(\mathsf X) \times \mathbb R$ be a point not contained in $\tilde{K}$. Suppose $W_1(\nu, \nu_0) \vee \lvert r - r_0 \rvert \geq \epsilon > 0, \forall (\nu, r) \in \tilde  K$. Then, there exists a tuple $(f,\alpha) \in [\mathscr L(\mathsf X)]_1 \times [-1, 1]$ satisfying $$\lVert f \rVert_{\lip} + \lvert \alpha \rvert \leq 1 \hspace{8pt}\text{and}\hspace{8pt}\int f d\nu + \alpha r  \geq \int f d\nu_0 + \alpha r_0 + \epsilon, \forall (\nu, r) \in \tilde K.$$
\end{theorem}
\proofbegin
 Let $\aes(\mathsf X) \times \mathbb R$ be equipped with the canonical norm $\lVert (m, r) \rVert := \lVert m \rVert_{\aes(\mathsf X)} \vee \lvert r \rvert.$  An extension of Theorem \ref{thm:aes} shows that its dual space is isometric to $\mathscr L(\mathsf X) \times \mathbb R$ with the norm $\lVert (f, r) \rVert = \lVert f \rVert_{\lip} + \lvert r \rvert,$ i.e. $(\aes(\mathsf X) \times \mathbb R)^* = \mathscr L(\mathsf X) \times \mathbb R$. The rest of the proof is similar to the proof of Theorem \ref{thm:mnd}.
\proofend

\subsection{A separation theorem on Wasserstein-1 space}

In this subsection, we extend Theorem \ref{thm:mnd} to the whole Wasserstein-1 space.
\begin{theorem}\label{thm:separation}
 Let $A$ be a convex and closed subset of $(\mathscr P_{1}(\mathsf X), W_1)$. Let $\mu_0 \in \mathscr P_{1}(\mathsf X)$ be a point not contained in $A$. Then, there exists a function $f \in [\mathscr L(\mathsf X)]_1$ such that
 $$
  \int f d \mu > \int f d\mu_0, \hspace{8pt}\forall \mu \in A.
 $$
\end{theorem}
\proofbegin

 Step 1. Since $A$ is closed and $\mu_{0}\notin A,$ there exists a positive  constant $\epsilon_0 > 0$ satisfying $W_1(\mu, \mu_0) \geq \epsilon_0, \forall \mu \in A.$

 Step 2. By Theorem 6.18 in \cite{villani2009optimal}, we can approximate $\mu$ by a point $\nu \in \mathscr D(\mathsf X)$ with any accuracy. Finally, Theorem \ref{thm:mnd} can be applied to find the required separation function in $\mathscr L(\mathsf X)$. More specifically, take $\epsilon_1 = \epsilon_0/5$ and define subsets of $\mathscr D(\mathsf X)$ as follows
 \begin{align*}
  B_{\epsilon_1}(\mu) := \{ \nu \in \mathscr D(\mathsf X) \mid W_1(\mu, \nu) < \epsilon_1 \},\hspace{8pt} A_{\epsilon_1} := \bigcup_{\mu \in A} B_{\epsilon_1}(\mu).
 \end{align*}
Since $\mathscr D(\mathsf X)$ is dense in $(\mathscr P_{1}(\mathsf X), W_1)$, $B_{\epsilon_1}(\mu) \neq \emptyset$, $\forall \mu \in \mathscr P_{1}(\mathsf X)$. Hence, for any $\nu_i \in A_{\epsilon_1}$ there exists a $\mu_i \in A$ such that $W_1(\nu_i, \mu_i) < \epsilon_1$, $i=1,2.$ Let $\nu_\alpha := \alpha \nu_1 + (1-\alpha) \nu_2$ and $\mu_\alpha := \alpha \mu_1 + (1-\alpha) \mu_2$, $\alpha \in (0,1)$. We have, by the dual representation \eqref{eq:kr-distance} of $W_1$, 
\begin{align*}
 W_1(\nu_\alpha, \mu_\alpha) \leq \alpha W_1(\nu_1, \mu_1) + (1-\alpha) W_1(\nu_2, \mu_2) < \epsilon_1,
\end{align*}
which implies that $A_{\epsilon_1}$ is convex. Similarly, we can find a $\nu_0 \in \mathscr D(\mathsf X)$ such that $W_1(\mu_0, \nu_0) < \epsilon_1$. Note that 
$
   W_1(\nu, \mu) \geq W_1(\mu_0, \mu) - W_1(\mu_0, \nu_0) - W_1(\nu_0, \nu)
                 \geq  \epsilon_0 - \epsilon_1 -  W_1(\nu_0, \nu), \forall \nu \in A_{\epsilon_1}, \mu \in A,
$
yields $\inf_{\mu \in A} W_1(\nu, \mu) \geq \epsilon_0 - \epsilon_1 -  W_1(\nu_0, \nu), \forall \nu \in A_{\epsilon_1}.$
Because $\epsilon_1 \geq \inf_{\mu \in A} W_1(\nu, \mu)$, we have $W_1(\nu_0, \nu) \geq \epsilon_0 - 2\epsilon_1 > 0, \forall \nu \in A_{\epsilon_1}.$ Hence, by Theorem \ref{thm:mnd}, there exists a function $f_0 \in \mathscr L(\mathsf X)$ with $ \lVert f_0 \rVert_{\lip} \leq 1,$ such that
\begin{align}
 \int f_0 d\nu \geq \int f_0 d\nu_0 + \epsilon_0 - 2\epsilon_1, \forall \nu \in A_{\epsilon_1}. \label{eq:1}
\end{align}

Step 3. Finally, we show that $f_0$ is the required Lipschitz function. For each $\mu \in A,$ there exists an $\nu_\mu \in A_{\epsilon_1}$ such that $W_1(\mu, \nu_\mu) < \epsilon_1$ and hence $\int f d \mu  \geq \int f d\nu_\mu -\epsilon_1, \forall f \in [\mathscr L(\mathsf X)]_1.$ In particular,  
\begin{align}
 \int f_0 d \mu \geq \int f_0 d \nu_\mu - \epsilon_1.\label{eq:2}
\end{align}
Similarly, $W_1(\mu_0, \nu_0) < \epsilon_1$ implies
\begin{align}
 \int f_0 d \nu_0 \geq \int f_0 d\mu_0 - \epsilon_1.\label{eq:3}
\end{align}
Combining \eqref{eq:1} -- \eqref{eq:3}, we obtain
$
 \int f_0 d \mu \geq \int f_0 d\mu_0 + \epsilon_0 - 4 \epsilon_1 = \int f_0 d\mu_0 + \frac{\epsilon_0}{5}, \forall \mu \in A,
$
which yields the required separability. \proofend

Analogously, we can obtain the same separation result in the space $\mathscr P_{1}(\mathsf X) \times \mathbb R$ which will be used in the following subsection. 
\begin{theorem}\label{thm:separation:extended}
 Let $\tilde A$ be a convex and closed subset of $\mathscr P_{1}(\mathsf X) \times \mathbb R$ equipped with the metric $$\tilde d((\mu, r_1),(\nu, r_2)) := W_1(\mu, \nu) \vee \lvert r_1 - r_2 \rvert.$$ Let $(\mu_0, r_0) \in \mathscr P_{1}(\mathsf X) \times \mathbb R$ be a point not contained in $A$. Then, there exists a tuple $(f, \alpha) \in \mathscr L(\mathsf X) \times \mathbb R,$ and $\epsilon_{0}>0,$ such that
 \begin{align*}
  \int f d \mu + \alpha r \geq \int f d\mu_0 + \alpha r_0 +\epsilon_{0}, \forall (\mu,r) \in \tilde A.
 \end{align*}
\end{theorem}

The proof is similar to the proof of Theorem \ref{thm:separation}. 

\subsection{Proof of the duality theorem}

We shall prove Theorem \ref{th:fm} in this subsection. We want to stress that the proofs of Lemma \ref{lm:M} and the Theorem \ref{th:fm} mostly follow the line of proof  found in \citet[ Theorem 2.2.15]{ioan2009duality}.

\begin{definition}\it
 For a function $\phi: \mathscr P_{1}(\mathsf X) \rightarrow \bar{\mathbb R}$, its \emph{domain} is defined by $\dom(\phi):= \{ \mu \in \mathscr P_{1}(\mathsf X): \phi(\mu) < +\infty \}$. $\phi$ is said to be \emph{proper} if $\phi(\mu) > - \infty$ for all $\mu \in \mathscr P_{1}(\mathsf X)$ and $\dom(\phi) \neq \emptyset$. The \emph{epigraph} of $\phi$ is defined as
$\epi(\phi) := \{ (\mu, r) \in \mathscr P_{1}(\mathsf X) \times \mathbb R \mid \phi(\mu) \leq r \}.$
\end{definition}

Now we provide the following proposition whose proof is standard and will be omitted.

\begin{proposition}\label{prop:epi}
 If $\phi: \mathscr P_{1}(\mathsf X) \rightarrow \bar{\mathbb R}$ is convex and lower semicontinuous (equipped with the metric $W_1$), then $\epi(\phi)$ is convex and closed w.r.t. $$\tilde d((\mu, a), (\nu, b)) := W_1(\mu, \nu) \vee \lvert a - b \rvert.$$
\end{proposition}

For any $\phi: \mathscr P_{1}(\mathsf X) \rightarrow \bar{\mathbb R},$ we define the following set
\begin{align}
M_{\phi}:=& \left\{ (f, \eta) \in \mathscr L(\mathsf X) \times \mathbb R \middle| \int f d\mu + \eta \leq \phi(\mu), \forall \mu \in \mathscr P_{1}(\mathsf X) \right\}.\label{eq:M}
\end{align}

\begin{lemma}\label{lm:M}
 Let $\phi: \mathscr P_{1}(\mathsf X) \rightarrow \bar{\mathbb R}$ be proper, convex and lower semicontinuous with respect to $W_1$. Then, $M_{\phi}$ is not empty.
\end{lemma}
\proofbegin
 Since $\phi$ is proper, there exists an element $\nu \in \mathscr P_{1}(\mathsf X)$ such that $\phi(\nu) \in \mathbb R$. Then $\epi(\phi) \neq \emptyset$ and $(\nu, \phi(\nu) - 1) \notin \epi(\phi).$ By Theorem \ref{thm:separation:extended}, there exists $(f_0, \eta_0) \in \mathscr L(\mathsf X) \times \mathbb R$ such that
 \begin{align*}
  \int f_0 d\nu + \eta_0(\phi(\nu) - 1) < \int f_0 d \mu + \eta_0 r, \forall (\mu, r) \in \epi(\phi).
 \end{align*}
 Since $(\nu, \phi(\nu)) \in \epi(\phi)$, we have $\eta_0 > 0$ and $(1/\eta_0) \left( \int f_0 d\nu - \int f_0 d \mu \right) + \phi(\nu) - 1 < r, \forall (\mu, r) \in \epi(\phi).$ For any $\mu \in \dom(\phi)$, we have $(\mu, \phi(\mu)) \in \epi(\phi),$ and hence $$(1/\eta_0)\left( \int f_0 d\nu - \int f_0 d \mu \right) + \phi(\nu) - 1 < \phi(\mu).$$ If  $\mu \notin \dom(\phi),$ this inequality holds trivially. Thus, $$\left(-f_0/\eta_0, \phi(\nu) - 1 + \int f_0 d\nu/\eta_0\right) \in M_{\phi}.$$\proofend 

It is easy to verify that under the same conditions as in the above lemma, we have the following equation for the second conjugate dual $\phi^c$ of $\phi$:
\begin{align*}
 \phi^c(\mu_0) = \sup_{(f,\eta) \in M_{\phi} }\left\{ \int f d\mu_0 + \eta \right\}.
\end{align*} 
We now proceed with the proof of Theorem \ref{th:fm}.

\begin{proof}[Proof of Theorem \ref{th:fm}]
 By definition, we have $\rho(f) \geq \int f d \mu - \phi(\mu), \forall (\mu, f) \in \mathscr P_{1}(\mathsf X) \times \mathscr L(\mathsf X),$ which implies
$\int f d \mu - \rho(f) \leq \phi(\mu), \forall (\mu, f) \in \mathscr P_{1}(\mathsf X) \times \mathscr L(\mathsf X),$ and hence $\phi^c(\mu) \leq \phi(\mu), \forall \mu \in \mathscr P_{1}(\mathsf X)$. 

Next we show that $\phi(\mu) \leq \phi^c(\mu), \forall \mu \in \mathscr P_{1}(\mathsf X)$. Assume towards a contradiction that there exists a $\mu_0 \in \mathscr P_{1}(\mathsf X)$ and $r_0 \in \mathbb R$ such that  
\begin{align}
 \phi(\mu_0) > r_0 > \phi^c(\mu_0) = \sup_{(f,\eta) \in  M_{\phi}}\left\{ \int f d\mu_0 + \eta \right\}.\label{eq:vcont}
\end{align}
It is clear that $(\mu_0, r_0) \notin \textrm{epi}(\phi)$. Note that by Proposition \ref{prop:epi}, $\epi(\phi)$ is closed. Furthermore it is convex and non-empty. Hence, by Theorem \ref{thm:separation:extended}, there exists a $(f_0, \eta_0)\in \mathscr L(\mathsf X) \times \mathbb R$ such that
\begin{align}
 \int f_0  d\mu + \eta_0 r  > \int f_0 d\mu_0 + \eta_0 r_0 + \epsilon, \forall (\mu, r) \in \epi(\phi), \label{eq:f_0}
\end{align}
Note that if $(\mu, r) \in \epi(\phi)$, then $(\mu, r+s) \in \epi(\phi)$ for any $s \geq 0$. Thus $\eta_0 \geq 0$. 

Suppose $\phi(\mu_0) \in \mathbb R$. Using the fact that $(\mu_0, \phi(\mu_0)) \in \epi(\phi)$, we obtain $\eta_0(\phi(\mu_0) - r_0) > 0$. Hence $\eta_0 > 0$. For $\mu \in \dom(\phi)$, we obtain
$\phi(\mu) > \frac{1}{\eta_0}\int f_0 d\mu_0 - \frac{1}{\eta_0} \int f_0 d\mu + r_0.$ Now setting $\eta = r_0 + \frac{1}{\eta_0}\int f_0 d\mu_0$ and $f = - \frac{f_0}{\eta_0}$ in \eqref{eq:vcont}, we have $\int f d\mu_0 + \eta = r_0,$  
which contradicts \eqref{eq:vcont}. 

Suppose now that $\phi(\mu_0) = +\infty.$ If $\eta_0 > 0$, the contradiction remains. Thus $\eta_0 = 0$ and \eqref{eq:f_0} becomes
$
 \int f_0  d\mu > \int f_0 d\mu_0 + \epsilon, \forall (\mu, r) \in \epi(\phi).
$
By Lemma \ref{lm:M}, the set $M_{\phi}$ defined in \eqref{eq:M} is not empty. Hence, there exists $(f_1, a_1) \in M_{\phi}$ such that $\int f_1 d\mu + a_1 \leq \phi(\mu), \forall \mu \in \mathscr P_{1}(\mathsf X).$ 
By the assumption made in \eqref{eq:vcont}, $b := (r_0 - \int f_1 d \mu_0 - a_1)/\epsilon > 0.$ For all $\mu \in \dom(\phi)$,
\begin{align}
  & \int (f_1 - b f_0) d\mu + \int  b f_0 d\mu_0 + a_1 + b \epsilon \label{eq:f_1} \\
 = & \int f_1 d \mu + a_1 + b \left( - \int f_0 d \mu + \int f_0 d \mu_0 + \epsilon \right) \leq \int f_1 d \mu + a_1 \leq \phi(\mu).\nonumber
\end{align}
This can be extended to all $\mu \in \mathscr P_{1}(\mathsf X)$. Hence, $(f_1 - b f_0, \int  b f_0 d\mu_0 + a_1 + b \epsilon) \in M_{\phi}.$ Taking $\mu = \mu_0$ in \eqref{eq:f_1}, we have $ \int (f_1 - b f_0) d\mu_0 + \int  b f_0 d\mu_0 + a_1 + b \epsilon = r_0,$ which again contradicts \eqref{eq:vcont}. 
\end{proof}

We will end this section with the proof Corollary \ref{coro:fm}

\begin{proof}[Proof of Corollary \ref{coro:fm}]
	(a) We show first $\phi(\mu) = \sup_{f \in \mathcal N_\phi} \int f d\mu$. Indeed, Theorem \ref{th:fm} yields 
	$$\phi(\mu) = \sup_{f \in \mathscr L(\mathsf X)} \left( \int f d \mu - \rho(f) \right) \geq \sup_{f \in \mathcal N_\phi} \left( \int f d \mu - \rho(f) \right) = \sup_{f \in \mathcal N_\phi} \int f d\mu.$$
	Since there exists at least one $\mu$ such that $\phi(\mu) \in \mathbb R$, we have $\rho(f) > -\infty, \forall f \in \mathscr L(\mathsf X)$. Thus,
	\begin{align*}
	\sup_{f \in \mathscr L(\mathsf X)} \left( \int f d \mu - \rho(f) \right) = \sup_{f \in \mathscr L(\mathsf X): \rho(f) < \infty} \left( \int f d \mu - \rho(f) \right) = \sup_{f \in \mathscr L(\mathsf X): \rho(f) \in \mathbb R} \left( \int f d \mu - \rho(f) \right).
	\end{align*}
	Due to the translation invariance, $f' := f -\rho(f)$ satisfies that $\rho(f') = 0$ if $\rho(f) \in \mathbb R$. Hence,
	\begin{align}
	\phi(\mu) = \sup_{f \in \mathscr L(\mathsf X): \rho(f) \in \mathbb R} \left( \int f d \mu - \rho(f) \right) = \sup_{f': f'=f-\rho(f), f\in \mathscr L(\mathsf X), \rho(f) \in \mathbb R} \int f' d \mu \leq \sup_{f' \in \mathcal N_\phi} \int f' d \mu. \label{eq:f'}
	\end{align}
	Combining the above two inequalities yields $\phi(\mu) = \sup_{f \in \mathcal N_\phi} \int f d\mu$.
	
	(b) We show now $\phi(\mu) = \sup_{f \in \bar{\mathcal N}_\phi} \int f d\mu$. Theorem \ref{th:fm} yields 
	$$\phi(\mu) = \sup_{f \in \mathscr L(\mathsf X)} \left( \int f d \mu - \rho(f) \right) \geq \sup_{f \in \bar{\mathcal N}_\phi} \left( \int f d \mu - \rho(f) \right) \geq \sup_{f \in \bar{\mathcal N}_\phi} \int f d\mu.$$
	On the other hand, \eqref{eq:f'} yields $\phi(\mu) \leq  \sup_{f' \in \mathcal N_\phi} \int f' d \mu \leq \sup_{f' \in \bar{\mathcal N}_\phi} \int f' d \mu.$ 
	\end{proof}

\section{POMDPs}
\subsection{Outline and related literature.} Most early literature on POMDPs (see e.g.\  \cite{sondik1978optimal, hernandez1989adaptive}) considers finite state spaces or bounded reward function. In their recent work, \citet{feinberg2016partially} consider a more general setting with Borel spaces and reward/cost functions that are bounded on one side. We adopt the setting of MDPs from \cite{hernandez1999further} with Borel spaces and real-valued reward functions, but we apply the weighted norm technique that connects them naturally to the Wasserstein-1 space.

In what follows, we are going to present the concept of POMDPs, and how they can be reduced to classical MDPs with a state space that contains all probability measures over the POMDPs' initial state space. We describe the weighted norm technique for solving MDPs, and demonstrate how our general assumptions imply the sufficient conditions for the convergence of the value iteration method in \cite{hernandez1999further}.

Although our approach covers cases with two-sided unbounded reward functions, these cases can be reduced to MDPs with bounded reward functions by applying an algebraic transformation  \citep{VanDerWal1981}. The contribution of our paper on POMDPs, is the conjugate approach for the computation of the  value function described in Section 5 that extends the ideas from \cite{sondik1978optimal, follmer2004stochastic}, and which we expect to inspire a computable set iteration algorithm for approximations of the value function.

\subsection{Setup}

\label{sec:pomdps}
A partially observable Markov decision process (POMDP, see e.g., \citet[Chapter 4]{hernandez1989adaptive}) is described by a tuple $(\mathsf X, \mathsf Y, \mathsf A, P, Q, \mu, r)$, where:
\begin{enumerate}[(a)]
 \item $\mathsf X$ is the (hidden or latent) state space, a Polish space with metric $d$. 
 \item $\mathsf Y$ is the space of observations,  a Borel space.
 \item $\mathsf A$ is the action space, a Borel space.
 \item $P(dx'|x, a)$ is the state transition law,  a stochastic kernel on $\mathsf X$ given $\mathsf K:= \mathsf X \times \mathsf A$. $\mathsf K$ is also a Borel space.
 \item $Q(dy |a, x)$ is the observation kernel a stochastic kernel on $\mathsf Y$ given $\mathsf K$. $Q_0$ is the initial observation kernel, a stochastic kernel on $\mathsf Y$ given $\mathsf \mathsf X$.
 \item $\mu \in \mathscr P(\mathsf X)$ is the initial distribution.
 \item $r: \mathsf K \rightarrow \mathbb R$ is the one-step reward function, which is $\mathcal B(\mathsf K)$-measurable.
\end{enumerate}

The POMDP evolves as follows; At time $t=0$, the initial (hidden or latent) state $x_0$ follows a given prior distribution $\mu$, while the initial observation $y_0$ is generated according to the initial observation kernel $Q_0(\cdot|x_0)$. If, at time $t,$ the state of the system is $x_t$ and the control $a_t \in \mathsf A$ is applied, then the agent receives a reward $r(x_t, a_t)$ and the system transits to state $x_{t+1}$ according to the transition law $P(dx_{t+1}| x_t, a_t)$. The observation $y_{t+1}$ is generated by the observation kernel $Q(dy_{t+1} | a_t, x_{t+1})$. The \emph{observed} history is defined as
\begin{align}
 h_0 := \{ \mu, y_0 \} \in \mathsf H_0 \textrm{ and } h_{t} := \{\mu, y_0, a_0, \ldots, y_{t-1}, a_{t-1}, y_t \} \in \mathsf H_t, t = 1,2, \ldots, \label{eq:h}
\end{align}
where $\mathsf H_0 := \mathscr P(\mathsf X) \times \mathsf Y$ and $\mathsf H_{t+1} = \mathsf H_t \times \mathsf Y \times \mathsf A,\, t=1,2,\ldots.$ Notably, comparing with the canonical Markov decision processes (MDPs, see e.g., \citet{hernandez1989adaptive}), the states $\{ x_t\}$ are not observable and hence, a \emph{policy} depends only on the observed history. 
%{\color{red} this is not true..still it might depend on the hidden states, I suggest just to erase {\it but is independent...}}{\color{blue} I have already modified accordingly.}

A deterministic policy $\boldsymbol \pi := [\pi_0, \pi_1, \ldots]$ is composed of a sequence of one-step policies $\pi_t: \mathsf H_t \rightarrow \mathsf A$, given the observed history up to time $t$. Let $\Pi$ be the set of all deterministic policies. Note that even with extension to nondeterministic policies, it is known \citep[Chapter 4]{hernandez1989adaptive} that an optimal policy to a POMDP is always deterministic. Hence, we consider in this paper only deterministic policies. The Ionescu-Tulcea theorem  \citep[pp.\ 140--141]{bertsekas1978stochastic} implies that for each $\pi \in \Pi$ and an initial $\mu \in \mathscr P(\mathsf X)$, along with $P$, $Q$ and $Q_0$, a probability measure $\mathbb P^{\boldsymbol \pi}_{\mu}$ and a stochastic process $\{X_t, Y_t, A_t\}$ can be defined in a canonical way. 
We denote by $\mathbb E^{\boldsymbol \pi}_{\mu}$ the expectation with respect to this probability measure $\mathbb P^{\boldsymbol \pi}_{\mu}$.

%defined by the products of Borel $\sigma$-algebras $\mathcal B(\mathsf X)$, $\mathcal B(\mathsf Y)$, and $\mathcal B(\mathsf A)$. %{\color{red} this should be written up more concrete... like eg: a  stochastic process $X_t, Y_t, A_t$ such that $P_\mu^\pi (X_{t+1} \in dx , Y_{t+1} \in dy , A_{t+1} ...| X_{0:t}, Y_{0:t}, A_{0:t}) = Q(dy |A_t, x) P (dx | X_t , A_t)$ ...and by the way isn't $A_t$ not $\pi_t$..?} {\color{blue} I have checked Feinberg's paper and revised the phrase slightly. I belief it should be okay now. In fact, the evolution of POMDPs have already been explained in the previous paragraph. Hence, the construction of $\mathbb P^{\boldsymbol \pi}_{\mu}$ should not be a problem...} 

We consider the following \emph{discounted cumulative rewards}
\begin{align}
 J_{T} (\boldsymbol \pi, \mu) := \mathbb E^{\boldsymbol \pi}_{\mu} \left[ \sum_{t=0}^T \alpha^t r(X_t, A_t) \right], \label{eq:t-stage}
\end{align}
where $\alpha \in (0,1)$ stands for a discount factor, and $T\in\mathbb{N}\cup\{\infty\}.$ 
The objective is now to maximize the expected reward over the set of deterministic policies $\Pi$,
\begin{align*}
 \phi_{T}^*(\mu) := \sup_{\pi \in \Pi} J_{T} (\boldsymbol \pi, \mu), \quad \mu \in \mathscr P(\mathsf X).
\end{align*}

We will finally use
 the notation \begin{equation}
\phi^{*}:=\phi^{*}_{\infty}.
\end{equation}

\subsection{Reduction to Markov decision process} We show briefly in this subsection that the POMDP can be reduced to a Markov decision process (MDP, see e.g.\ \citet{hernandez1999further}). We follow mostly the derivation by \citet[Chapter 4]{hernandez1989adaptive}. We first introduce the following notation:
\begin{align}
 \tilde r(\mu, a) := \int r(x, a) \mu(dx), \ \textrm{ and } \ \tilde P(\mathsf B|\mu,a) := 
 \int  \mu(dx) P(\mathsf B|x,a), \mathsf B \in \mathcal B(\mathsf X),
\end{align}
where $\mu\in \mathscr P(\mathsf X),$ and $a\in \mathsf A.$
For any $\mathsf C \in \mathcal B(\mathsf Y)$ and $\mathsf B \in \mathcal B(\mathsf X)$, we define
\begin{align}
 R(\mathsf B, \mathsf C | \mu, a) := \int_{\mathsf B} Q(\mathsf C | a, x') \tilde P(dx'|\mu, a), \ \textrm{ and } \ 
 \tilde R(\mathsf C | \mu, a) := &  R(\mathsf X, \mathsf C | \mu, a).
\end{align}

\begin{proposition}
 There exists a stochastic kernel $M$ from $\mathscr P(\mathsf X) \times \mathsf A \times \mathsf Y$ to  $\mathsf X$ such that for each $\mu \in \mathscr P(\mathsf X)$, $a \in \mathsf A, \mathsf B \in \mathcal B(\mathsf X)$ and $\mathsf C \in \mathcal B(\mathsf Y)$,
 \begin{align}
  R(\mathsf B, \mathsf C | \mu, a) = \int_{\mathsf C} M(\mathsf B | \mu, a, y) \tilde R(dy | \mu, a). \label{eq:RM}
 \end{align}
\end{proposition}
\proofbegin 
Direct application of \citet[Corollary 7.27.1]{bertsekas1978stochastic}. 
\proofend

%\end{proof}

$M$ can also be viewed as a mapping $\mathscr P(\mathsf X) \times \mathsf A \times \mathsf Y \rightarrow \mathscr P(\mathsf X)$. With a slight abuse of notation, let $M(\mu, a, y) := M(\cdot | \mu, a, y) \in \mathscr P(\mathsf X).$ Then we define the following stochastic kernel
\begin{align}
 \tilde Q(\mathsf D | \mu, a) := \int \mathsf 1_{\mathsf D}(M(\mu, a, y)) \tilde R(dy|\mu, a),\hspace{8pt} D \in \mathcal{B}(\mathscr P(\mathsf X)). \label{eq:tildeQ}
\end{align}

Let $\mu_t \in \mathscr P(\mathsf X)$ be the distribution at time $t.$ Then, given an action $a_t \in \mathsf A$ and an observation $y_{t+1} \in \mathsf Y$, the \emph{successive} distribution $\mu_{t+1} \in \mathscr P(\mathsf X)$ is given by
\begin{align}
 \mu_{t+1} = M(\mu_t, a_t, y_{t+1}). \label{eq:Mmu}
\end{align}
Note that $\mu_{t+1}$ is a random measure since $y_{t+1}$ is a random variable with the distribution $\tilde R(\cdot|\mu_t, a_t)$. Hence, the POMDP can be reduced to a Markov decision process (MDP) with the \emph{belief} state space $\mathscr P(\mathsf X)$, the action space $\mathsf A$, the reward function $\tilde r$ on $\mathscr P(\mathsf X) \times \mathsf A$ and the transition kernel on belief states $\tilde Q$ defined above. %For a formal definition on MDPs on Borel spaces with weighted norm, we refer to \cite{hernandez1999further} and Appendix \ref{sec:mdps}. %Note that the metric for $\mathscr P$ is not specified yet, which is postponed to the next subsection.

Let $\tilde h_t$ be a $t$-stage history for the MDP described above:
$$
  \tilde h_t := \{   \tilde\mu_0,   \tilde a_0, \ldots,   \tilde \mu_{t-1},   \tilde a_{t-1},   \tilde \mu_t \} \in \tilde{\mathsf H}_t := \mathscr P(\mathsf X) \times (\mathsf A \times \mathsf P)^t, t = 0, 1, \ldots,
$$
where $  \tilde\mu_0$ is the initial distribution and the $\tilde{\mu}_t \in \mathscr P(\mathsf X)$ are recursively defined by \eqref{eq:Mmu}. Given the original $t$-stage history $h_t$ defined in \eqref{eq:h}, let $m_t: \mathsf H_t \rightarrow \tilde{\mathsf H}_t$ be the mapping such that $m_t(h_t) = \tilde h_t, \forall h_t \in \mathsf H_t, t= 0, 1, \ldots.$ For a history-dependent MDP-policy $\boldsymbol \delta := [\delta_0, \delta_1, \ldots]$, where $\delta_t: \tilde{\mathsf H}_t \rightarrow \mathsf A$, we define its counterpart POMDP-policy as $\boldsymbol \pi^\delta = [\pi_0^\delta, \pi_1^\delta, \ldots]$, with
$
\pi_t^\delta(h_t) = \delta_t(m_t(h_t)).
$
If a policy $\boldsymbol \delta$ is optimal for the MDP, then its counterpart policy $\boldsymbol \pi^\delta$ is also optimal for the POMDP. For more details, we refer to \citet[Chapter 4]{hernandez1989adaptive} and references therein.

After the reduction to a MDP, it is well known that under proper assumptions (see, e.g., \cite[Chapter 4]{hernandez1989adaptive} and \cite{feinberg2016partially}), the optimal $\phi^*$ for the infinite-stage case satisfy the following \emph{optimality equation}:
\begin{align}
 \phi(\mu) = \mathcal T(\phi)(\mu) := \sup_{a \in \mathsf A} \left( r(\mu, a) + \alpha \int \phi\left(M(\mu, a, y) \right) \tilde R(dy|\mu, a) \right), \forall \mu \in \mathscr P(\mathsf X), \label{eq:optimal}
\end{align}
where $\mathcal T$ is an operator on the space of Borel measurable functions on $\mathscr P(\mathsf X)$. %We define
%\begin{align*}
% \mathcal T_a(\phi)(\mu) := r(\mu, a) + \alpha \int \phi\left(M(\mu, a, y) \right) \tilde R(dy|\mu, a).
%\end{align*}
%Then, we have $\mathcal T(\phi)(\mu) = \sup_{a \in \mathsf A} \mathcal T_a(\phi)(\mu), \forall \mu \in \mathscr P.$ 
We will show in Section \ref{sec:optimal}, that under general assumptions, the existence of a solution of the above equation is guaranteed, as well as the existence of an optimal deterministic policy. The value iteration algorithm gives a sequence that converges to the value function.

\section{Optimal solution on weighted space}\label{sec:optimal}
\subsection{Weighted norm}

The \emph{weighted norm} has proven to be very useful when dealing with MDPs (see e.g.\ \cite{hernandez1999further}). The weighted norm induced by the Wasserstein-1 metric is defined as follows: First, we specify a weight function $w: \mathsf X \rightarrow [1, \infty)$ by
\begin{align}
w(x) := 1 + k \cdot d(x_0, x),\hspace{4pt} x \in \mathsf X \label{eq:weightfunc}
\end{align}
with some fixed $x_0 \in \mathsf X$ and a positive constant $k > 0$. This definition implies $w(x_0) = 1$. One can easily verify that $w$ is Lipschitz, and therefore continuous and measurable. Denote by $\mathscr L_{w}(\mathsf X)$ the space of all continuous functions on $\mathsf X$ such that
\begin{align*}
\lVert f \rVert_w := \sup_{x \in \mathsf X} \frac{\lvert f(x) \rvert}{w(x)} < \infty. 
\end{align*}
Let $\mathscr P_w (\mathsf X)\subset \mathscr P(\mathsf X)$ be the set of probability measures  $\mu$ on $\mathsf X$ satisfying $\int w d \mu < \infty$. We equip $\mathscr P_w (\mathsf X)$ with the following concept of weak convergence: 

\begin{definition}\label{def:weak}\it
	(a) $\mu_n$ is said to converge weakly to $\mu$ in $\mathscr P_w(\mathsf X)$, if (i) for any bounded continuous function $f$ on $\mathsf X$, it follows that $ \int_{\mathsf X} f(x) \mu_n(dx) \rightarrow \int_{\mathsf X} f(x) \mu(dx) \textrm{ as } n \rightarrow \infty;$ 
	and (ii) $\int_{\mathsf X} w(x) \mu_n(dx) \rightarrow \int_{\mathsf X} w(x) \mu(dx)$ as $n \rightarrow \infty$. (b) A function $\phi: \mathscr P_w(\mathsf X)\rightarrow \mathbb R$ is said to be lower (resp.\ upper) semicontinuous if $$\liminf_{n\rightarrow \infty} \phi(\mu_n) \geq \phi(\mu) \ (resp.\ \limsup_{n\rightarrow \infty} \phi(\mu_n) \leq \phi(\mu)),$$ whenever $\mu_n$ converges weakly to $\mu$. $\phi$ is said to be continuous, if $\phi$ is both lower and upper continuous.
\end{definition}

\begin{remark}
	It is worth to be mentioned that the weak convergence defined above is stronger than the usual weak convergence, which requires only (i). Hence, to emphasize this difference, we call the latter 
	\emph{canonical weak convergence} throughout the rest of this paper.
\end{remark}

%It is remarkable that the definition of $w$-weak convergence is different from the usual one. The latter requires only condition (i). 

\begin{proposition}\label{prop:weakconv}
	The following two statements are equivalent: (i) $\mu_n$ converges weakly to $\mu$ in $\mathscr P_w(\mathsf X)$ and (ii) $\int f d\mu_n \rightarrow \int f d \mu, \forall f \in \mathscr L_w(\mathsf X)$.
\end{proposition}
For a proof, see \citet[Definition 6.8]{villani2009optimal}. 

\begin{proposition}
	Let $w$ be defined as in \eqref{eq:weightfunc} with some $k > 0$ and $x_0 \in \mathsf X$. Then (i) $\mathscr P_{1}(\mathsf X) = \mathscr P_w(\mathsf X)$ and (ii) the weak convergence defined in Definition \ref{def:weak} is equivalent to the convergence in $(\mathscr P_{1}(\mathsf X), W_1)$, in other words, $W_1$ metrizes the weak convergence.
\end{proposition}
\proofbegin (i) This is obvious by definition. (ii) This is a direct result of \citet[Theorem 6.9]{villani2009optimal}. \proofend

In the rest of this paper, we always assume that the weight function satisfies \eqref{eq:weightfunc} and hence $\mathscr P_{w}(\mathsf X)$ is used interchangeably with $\mathscr P_{1}(\mathsf X)$. Then, the belief state space is Polish, and therefore, a Borel space. 

We next specify the weight function and its weighted norm on $\mathscr P_w(\mathsf X)$. Define $\tilde w : \mathscr P_w(\mathsf X) \rightarrow [1, \infty)$ as
\begin{align}
 \tilde w(\mu) := \int w d\mu. \label{eq:tildew}
\end{align}
It is easy to check that $\tilde w$ is a continuous function and hence measurable on $(\mathscr P_w(\mathsf X), W_1)$. Define the following space of functions on $\mathscr P_w(\mathsf X)$ with bounded $\tilde w$-norm:
\begin{align*}
 \mathscr B_{\tilde w}(\mathsf X) := \left\{ \phi : \mathscr P_w(\mathsf X) \rightarrow \mathbb R \ \middle| \ \phi \textrm{ is } \mathcal B(\mathscr P_w(\mathsf X))\textrm{-measurable}, \lVert \phi \rVert_{\tilde w} := \sup_{\mu \in \mathscr P_w(\mathsf X)}\frac{\lvert \phi(\mu) \rvert}{\tilde w(\mu)} < \infty\right\}. 
\end{align*}
In the next subsection, we shall specify some assumptions on the original POMDP in order to ensure the  assumptions needed for MDPs as in  Theorem 8.3.6 in \citet{hernandez1999further}.

\subsection{Assumptions} We introduce the following assumption for the reward function.
\begin{assumption}\label{asmp:reward} (i) There exists a positive constant $\bar r > 0$ such that $\lvert r(x,a) \rvert \leq \bar r w(x)$, for each $(x,a) \in \mathsf K.$ (ii) For each $x \in \mathsf X$, $a \mapsto r(x,a)$ is upper semicontinuous.
\end{assumption}

\begin{proposition}\label{prop:reward}
 Under Assumption \ref{asmp:reward}, (i) $\lvert \tilde{r}(\mu,a) \rvert \leq \bar r \tilde w(\mu), \forall (\mu,a) \in \mathscr P_w(\mathsf X) \times \mathsf A$; (ii) for each $\mu \in \mathscr P_w(\mathsf X)$, $a \mapsto \tilde r(\mu, a)$ is upper semicontinuous; (iii) for each $a \in \mathsf A$, $\mu \mapsto \tilde r(\mu, a)$ is continuous in $\mathscr P_w(\mathsf X)$. 
\end{proposition}
\proofbegin
(i) For each $(\mu,a) \in \mathscr P_w(\mathsf X) \times \mathsf A$, we have
$\lvert \tilde{r}(\mu,a) \rvert \leq \int \lvert r(x,a) \rvert \mu(dx) \leq \bar r \int w d\mu = \bar r \tilde w(\mu).$ (ii) Let $\{ a_n, n=1, 2,\ldots \}$ be a sequence of actions converging to $a_0$ and set $r_n(x) := r(x,a_n), n \in \mathbb N$. By Assumption \ref{asmp:reward}(i), $r_n \leq \bar r w$. Applying the reversed Fatou's lemma, we obtain
 \begin{align*}
  \limsup_{n \rightarrow \infty} \int r_n d\mu \leq \int \limsup_{n \rightarrow \infty} r_n d\mu \leq \int r_0 d \mu,
 \end{align*}
where the last inequality is due to Assumption \ref{asmp:reward}(ii). Finally, (iii) is a direct result of Proposition \ref{prop:weakconv}(ii). \proofend

Similar to the assumptions made in the literature of MDPs \citep[Assumptions 8.3.2 and 8.3.3]{hernandez1999further}, we introduce the following assumption on the transition kernel $P$:
\begin{assumption}\label{asmp:P}
 \begin{enumerate}[(i)]
  \item There exists a constant $\beta \in (0, \alpha^{-1})$ such that $$\int w(x') P(dx'|x,a) \leq \beta w(x), \forall (x,a) \in \mathsf X \times \mathsf A.$$
  \item For each $x \in \mathsf X$, $a \mapsto \int w(x')P(dx'|x,a)$ is continuous.
 \end{enumerate}
\end{assumption}

Under the above assumption, we show that the new probability measure $M(\mu, a, y)$ belongs to $\mathscr P_w(\mathsf X)$ almost surely, and $a\mapsto\int \tilde w(\mu') \tilde Q(d\mu'|\mu,a)$ is continuous. 

\begin{proposition}\label{prop:transition}
 Suppose Assumption \ref{asmp:P} holds. Then for each $\mu \in \mathscr P_w(\mathsf X)$ and $a \in \mathsf A$, (i) $M(\mu, a, y) \in \mathscr P_w(\mathsf X)$, $\tilde R(\cdot|\mu, a)$-almost surely; (ii) $\int \tilde w(\mu') \tilde Q(d\mu'|\mu,a) \leq \beta \tilde w(\mu);$ and (iii) for each $\mu \in \mathscr P_w(\mathsf X)$, the map $a \mapsto \int \tilde w(\mu') \tilde Q(d\mu'|\mu,a)$ is continuous.
\end{proposition}
\proofbegin Fix an arbitrary $(\mu, a) \in \mathscr P_w (\mathsf X)\times \mathsf A$. (i) Let $\mathsf C \in \mathcal B(\mathsf Y)$ be a subset such that $\tilde R(\mathsf C|\mu, a) > 0$. Then, we have
 \begin{align*}
  \int_{\mathsf C} \int_{\mathsf X} w(x') M(dx'|\mu, a, y) \tilde R(dy|\mu, a) = \int_{\mathsf X} Q(\mathsf C | a, x') w(x') \tilde P(dx'|\mu, a)  \leq \\  \int_{\mathsf X} Q(\mathsf Y | a, x') w(x') \tilde P(dx'|\mu, a)
  \leq \int_{\mathsf X} w(x') \tilde P(dx'|\mu, a) \leq \beta \int w d \mu < \infty.
 \end{align*}
This implies that $\int_{\mathsf X} w(x') M(dx'|\mu, a, y) < \infty$, $\tilde R(\cdot|\mu, a)$-almost surely, and hence (i) holds.

(ii) By definition, we have
 \begin{equation*}
\begin{split}	
	\int \tilde w(\mu') \tilde Q(d\mu'|\mu,a) = & \int_{\mathsf Y} \tilde w(M(\mu, a, y)) \tilde R(dy | \mu, a) \\
	\textrm{(by \eqref{eq:tildew})}\quad = & \int_{\mathsf Y} \int_{\mathsf X} w(x')M(dx'| \mu, a, y) \tilde R(dy | \mu, a) \\
	\textrm{(by \eqref{eq:RM} and Fubini's theorem)}\quad  = & \int_{\mathsf X} Q(\mathsf Y | a, x') w(x') \tilde P(dx'|\mu, a) \\
	\leq & \beta \int w d\mu = \beta \tilde w(\mu).
\end{split}
\end{equation*}

 (iii) Note that the above calculation yields
 $$\int_{\mathsf X} \tilde w(\mu') \tilde Q(d\mu'|\mu,a) = \int_{\mathsf X} w(x') \tilde P(dx'|\mu, a) = \int_{\mathsf X}\int_{\mathsf X} w(x') P(dx'|x, a) \mu(dx).$$
 Let $(a_n)_{n\in\mathbb{N}}$ be a sequence in $\mathsf{A}$ converging to $a_0$ and define $$f_n(x) := \int w(x') P(dx'|x, a_n), n \in \mathbb N.$$ Hence, the required continuity is equivalent to showing that $\lim_{n \rightarrow \infty}\int f_n d\mu = \int f_0 d\mu.$ Indeed, by Assumption \ref{asmp:P}(i), we have $f_n \leq \beta w, \forall n \in \mathbb N$. The reversed Fatou's lemma implies $\limsup_{n \rightarrow \infty} \int f_n d\mu \leq \int \limsup_{n \rightarrow \infty} f_n d\mu = \int f_0 d\mu.$ 
 On the other hand, we have $f_n \geq -\beta w, \forall n \in \mathbb N$. Then, the extended Fatou's lemma implies $
  \liminf_{n \rightarrow \infty} \int f_n d\mu \geq \int \liminf_{n \rightarrow \infty} f_n d\mu = \int f_0 d\mu.$
 Combining the above two inequalities yields the convergence. \proofend

\begin{assumption}\label{asmp:PQ} For each $\mu \in \mathscr P_w(\mathsf X)$, there exist stochastic kernels $M$ on $\mathscr P(\mathsf X)$ given $\mathscr P_w(\mathsf X) \times \mathsf A \times \mathsf Y$ and $\tilde R$ on $\mathsf Y$ given $\mathscr P_w(\mathsf X) \times \mathsf A$ satisfying \eqref{eq:RM} such that, if $\{a_n \in \mathsf A, n=1,2,\ldots\}$ converges to $a_0 \in \mathsf A$ as $n \rightarrow \infty$, 
 \begin{enumerate}[(i)]
  \item there exists a subsequence $\{ a_{n_k} \} \subset \{ a_{n} \}$ and a measurable set $\bar{\mathsf C} \in \mathcal B(\mathsf Y)$ such that $\tilde R(\bar{\mathsf C} | \mu, a_0) = 1$ and for all $y \in \bar{\mathsf C}$, $M(\mu, a_{n_k}, y)$ converges canonically weakly to $M(\mu, a_0, y)$;
  \item for each $\mathsf C \in \mathcal B(\mathsf Y)$, $\tilde R(\mathsf C  |\mu, a_n) \rightarrow \tilde R(\mathsf C  |\mu, a_0)$ as $n \rightarrow \infty$.
 \end{enumerate}
\end{assumption}
This assumption is inspired by Condition (c) in \cite[Theorem 3.2]{feinberg2016partially}. A sufficient condition for it will be discussed in the next section (see Remark \ref{rm:ascont}).

\begin{proposition}\label{prop:tildeQ}
 Under Assumptions \ref{asmp:PQ} and \ref{asmp:P}, for each $\mu \in \mathscr P_w(\mathsf X)$, $a \mapsto \tilde Q(\cdot | \mu, a)$ is canonically weakly continuous. %, where $\tilde Q$ and $\tilde w$ are defined in \eqref{eq:tildeQ} and \eqref{eq:tildew} respectively.
\end{proposition}
\proofbegin This canonical weak continuity is guaranteed by \citet[Theorem 3.4]{feinberg2016partially}. %It remains to show that $a \mapsto \int \tilde w(\mu') \tilde Q(d\mu'|\mu, a)$ is continuous, which is guaranteed by Proposition \ref{prop:transition}(iii). 
\proofend

Finally, to guarantee the existence of one ``selector'', we assume
\begin{assumption}\label{asmp:A}
 $\mathsf A$ is compact.
\end{assumption}
Note that the operator $\mathcal T: \mathscr B_{\tilde w}(\mathsf X) \rightarrow \mathscr B_{\tilde w}(\mathsf X)$ is defined as follows
\begin{align}
 \mathcal T_a(\phi)(\mu) := & \tilde r(\mu, a) + \alpha \int \phi\left(M(\mu, a, y) \right) \tilde R(dy|\mu, a) \label{eq:ta}\\
 \quad \textrm{and} \quad \mathcal T(\phi)(\mu) := & \sup_{a \in \mathsf A} \mathcal T_a(\phi)(\mu). \label{eq:t}
\end{align}
Under Assumptions \ref{asmp:reward} -- \ref{asmp:PQ}, it is guaranteed that for each $\mu \in \mathscr P_w(\mathsf X)$, $a \mapsto \mathcal T_a(\phi)(\mu)$ is upper-semicontinuous (for a proof, see \cite[Lemma 8.3.7(a)]{hernandez1999further}). Hence, under the additional Assumption \ref{asmp:A}, the optimal $a$ in the above optimization problem is always attainable in $\mathsf A$ (see, e.g., \cite[Lemma 8.3.8(a)]{hernandez1999further}). Hence, from now on, we replace ``sup'' with ``max''.

\subsection{Value iteration} 

The following \emph{value iteration} is a widely used method to compute the optimal solution for POMDPs, and MDPs as well. Starting from arbitrary \emph{value function} in $\mathscr B_{\tilde w}(\mathsf X)$, $\phi_0$, at time $t$, we update value function as follows
\begin{align*}
 \phi_{t+1}(\mu) = \mathcal T (\phi_t)(\mu) =  \max_{a \in \mathsf A} \left( \tilde r(\mu, a) + \alpha \int \phi\left(M(\mu, a, y) \right) \tilde R(dy|\mu, a)  \right) 
\end{align*}
Finally, by a suitable application of Theorem 8.3.6 in \citet{hernandez1999further}, we obtain the following convergence.
\begin{theorem}\label{thm:vi}
 Suppose that Assumptions \ref{asmp:reward} -- \ref{asmp:A} hold. Let $\beta$ be the constant in Assumption \ref{asmp:P}(i) and $\bar r$ be the constant in \ref{asmp:reward}(i) and define $\gamma := \alpha \beta \in (0,1)$. Then
 \begin{enumerate}[(a)]
  \item the optimal value function $\phi^*$ is the unique fixed point of the operator $\mathcal T$ satisfying $\phi^* = \mathcal T(\phi^*)$ in $\mathscr B_{\tilde w}(\mathsf X)$ and 
  $
   \lVert \phi_t - \phi^*\rVert_{\tilde w} \leq \bar r \gamma^t/(1-\gamma), t = 1, 2, \ldots.
  $
  \item there exists a selector $f^*: \mathscr P_w(\mathsf X) \rightarrow \mathsf A$ such that 
  $$
  \phi^*(\mu) = \tilde r(\mu, f^*(\mu)) + \alpha \int \phi\left(M(\mu, f^*(\mu), y) \right) \tilde R(dy|\mu, f^*(\mu)), \forall \mu \in \mathscr P_w(\mathsf X).
  $$ and $\boldsymbol \pi^* = (f^*)^\infty$ is one optimal policy satisfying $\phi^*(\mu) = J(\mu, \boldsymbol \pi^*), \forall \mu \in \mathscr P_w(\mathsf X)$.
 \end{enumerate}
\end{theorem}
\proofbegin The original POMDP specified in Subsection \ref{sec:pomdps} can be reduced to an MDP with $(\mathscr P_w(\mathsf X), \mathsf A, \tilde r, \tilde Q)$. %defined in Appendix \ref{sec:mdps}. 
Under Assumptions \ref{asmp:reward} -- \ref{asmp:A}, Propositions \ref{prop:reward} -- \ref{prop:tildeQ} hold, and therefore, the conditions required by Theorem 8.3.6 in \citet{hernandez1999further} are satisfied. The assertion is then a direct application of that Theorem.\proofend

\subsection{Application to POMDPs}

Now we apply the conjugate duality obtained in Corollary \ref{coro:fm} to POMDPs. Recall that the operators $\mathcal T_a$ and $\mathcal T$ are defined in equations \eqref{eq:ta} and \eqref{eq:t}.
\begin{lemma}\label{lm:convex}
 If $\phi: \mathscr P_w(\mathsf X) \rightarrow \mathbb R$ is convex, then $\mathcal T_a (\phi)$ is convex, $\forall a \in \mathsf A$, and therefore, $\mathcal T(\phi)$ is convex as well.
\end{lemma}
\proofbegin
It is sufficient to show that $\mu \mapsto \tilde{r}(\mu,a) + \alpha \int \phi(M(\mu, a, y)) \tilde R(dy|\mu, a)$ is convex for each $a \in \mathsf A$. Indeed, take any action $a\in \mathsf A$ and let $\mu_1$ and $\mu_2$ be two arbitrary elements in $\mathscr P_w(\mathsf X)$. Take any $\kappa \in (0,1)$ and define $\mu_\kappa := \kappa \mu_1 + (1-\kappa) \mu_2$. By the definition, for any $\mathsf B \in \mathcal B(\mathsf X)$ and $\mathsf C \in \mathcal B(\mathsf Y)$, we have
 \begin{align}
  R(\mathsf B, \mathsf C | \mu_\kappa, a) = & \kappa R(\mathsf B, \mathsf C | \mu_1, a) + (1-\kappa) R(\mathsf B, \mathsf C | \mu_2, a)  \\
                                   = & \kappa \int_{\mathsf C} M(\mathsf B | \mu_1, a, y) \tilde R(dy | \mu_1, a) + (1-\kappa) \int_{\mathsf C} M(\mathsf B | \mu_2, a, y) \tilde R(dy | \mu_2, a). \label{eq:rbc}
 \end{align}
On the other hand, a simple calculation yields
 \begin{align}
  \tilde R(\mathsf C |\mu_\kappa, a) = \kappa \tilde R(\mathsf C |\mu_1, a) + (1-\kappa) \tilde R(\mathsf C |\mu_2, a), \forall \mathsf C \in \mathcal B(\mathsf Y). \label{eq:rkappa}
 \end{align}
 Hence, $\tilde R(\mathsf C |\mu_\kappa, a) = 0$ implies $\tilde R(\mathsf C |\mu_1, a) = 0$ and $\tilde R(\mathsf C |\mu_2, a) = 0, \forall \mathsf C \in \mathcal B(\mathsf Y)$. By Radon-Nikodym theorem, there exist functions $f_i: \mathsf Y \times \mathsf A \rightarrow [0, \infty), i = 1, 2$, which are both $\mathcal B(\mathsf Y)$-measurable for the fixed $a$, such that 
 \begin{align}
  \tilde R(\mathsf C |\mu_i, a) =  \int_{\mathsf C} f_i(y, a) \tilde R(dy |\mu_\kappa, a), i = 1, 2. \label{eq:rtildec}
 \end{align}
 Applying these two equations in \eqref{eq:rkappa} accordingly, we obtain $$\tilde R(\mathsf C |\mu_\kappa, a) = \int_{\mathsf C} \left( \kappa f_1(y,a) + (1-\kappa) f_2(y,a) \right)  \tilde R(dy |\mu_\kappa, a),  \forall \mathsf C \in \mathcal B(\mathsf Y),$$ which implies that $\kappa f_1(y,a) + (1-\kappa) f_2(y,a) = 1$, $\tilde R(\cdot |\mu_\kappa, a)$-almost surely. In other words, there exists a Borel set $\bar{\mathsf C} \in \mathcal B(\mathsf Y)$ such that $\tilde R(\bar{\mathsf C}|\mu_\kappa, a) = 1$ and $\kappa f_1(y,a) + (1-\kappa) f_2(y,a) = 1, \forall y \in \bar{\mathsf C}.$
 
Applying \eqref{eq:rtildec} to \eqref{eq:rbc}, we obtain
\begin{align*}
  R(\mathsf B, \mathsf C | \mu_\kappa, a) = & \int_{\mathsf C} \big[ \kappa M(\mathsf B | \mu_1, a, y) f_1(y, a) + (1-\kappa) M(\mathsf B | \mu_2, a, y) f_2(y, a) \big] \tilde R(dy |\mu_\kappa, a),
\end{align*}
and for each $\kappa \in (0,1)$, $M(\cdot | a, y, \kappa) := \kappa M(\cdot | \mu_1, a, y) f_1(y, a) + (1-\kappa) M(\cdot | \mu_2, a, y) f_2(y, a)$ is a valid stochastic kernel satisfying
$R(\mathsf B, \mathsf C | \mu_\kappa, a) = \int_{\mathsf C} M(\mathsf B | a, y, \kappa) \tilde R(dy |\mu_\kappa, a), \forall \mathsf B \in \mathcal B(\mathsf X), \mathsf C \in \mathcal B(\mathsf Y).$
Finally, the convexity of $\phi$ implies
\begin{align*}
  & \int \phi(M(\cdot|a,y, \kappa)) \tilde R(dy |\mu_\kappa, a) = \int_{\bar{\mathsf C}} \phi(M(\cdot|a,y, \kappa)) \tilde R(dy |\mu_\kappa, a)\\
 \leq & \int_{\bar{\mathsf C}} \left[ \kappa f_1(y, a) \phi( M(\mu_1, a, y) ) + (1-\kappa) f_2(y, a) \phi( M( \mu_2, a, y) ) \right]\tilde R(dy |\mu_\kappa, a) \\
 \leq & \kappa \int \phi( M( \mu_1, a, y) ) \tilde R(dy |\mu_1, a) + (1-\kappa) \int \phi( M( \mu_2, a, y) ) \tilde R(dy |\mu_2, a),
\end{align*}
which yields the required convexity.
\proofend

We introduce the following assumption accompanying Assumption \ref{asmp:PQ}. 
\begin{assumption}\label{asmp:mucontinuity} For each $a \in \mathsf A$, there exist stochastic kernels $M$ on $\mathscr P(\mathsf X)$ given $\mathscr P_w(\mathsf X) \times \mathsf A \times \mathsf Y$ and $\tilde R$ on $\mathsf Y$ given $\mathscr P_w (\mathsf X)\times \mathsf A$ satisfying \eqref{eq:RM} such that, if $\{\mu_n \in \mathscr P_w(\mathsf X), n=1,2,\ldots\}$ converges to $\mu_0 \in \mathscr P_w(\mathsf X)$ as $n \rightarrow \infty$, 
 \begin{enumerate}[(i)]
  \item there exists a subsequence $\{ \mu_{n_k} \} \subset \{ \mu_{n} \}$ and a measurable set $\bar{\mathsf C} \in \mathcal B(\mathsf Y)$ such that $\tilde R(\bar{\mathsf C} | \mu_0, a) = 1$ and for all $y \in \bar{\mathsf C}$, $M(\mu_{n_k}, a, y)$ converges canonically weakly to $M(\mu_0, a, y)$;
  \item for each $\mathsf C \in \mathcal B(\mathsf Y)$, $\tilde R(\mathsf C  |\mu_n, a) \rightarrow \tilde R(\mathsf C  |\mu_0, a)$ as $n \rightarrow \infty$.
 \end{enumerate}
\end{assumption}

\begin{lemma}\label{lm:lsc}
Suppose Assumptions \ref{asmp:mucontinuity}, \ref{asmp:reward}(i) and \ref{asmp:P}(i) hold. %Assume futher that there exists a constant $\bar \phi > 0$ such that $ \lvert \phi(\mu) \rvert \leq \bar \phi \tilde w(\mu) = \bar \phi \int w d\mu.$ 
Then, for any $\phi \in \mathscr B_{\tilde w}(\mathsf X)$, $\mathcal T_a (\phi)$ is continuous for any $a \in \mathsf A$ and therefore $\mathcal T(\phi)$ is lower semicontinuous.
\end{lemma}
\proofbegin
Fix one $a \in \mathsf A$. Note that Proposition \ref{prop:reward}(iii) ensures the continuity of the function $\mu \mapsto r(\mu, a)$. It remains to show that $\mu \mapsto \int \phi(\mu') \tilde Q(d\mu' | \mu, a)$ is continuous.

By Assumption \ref{asmp:mucontinuity} and \cite[Theorem 3.4]{feinberg2016partially}, $\mu \mapsto \tilde Q(d \mu'|\mu, a)$ is canonically weakly continuous. By Assumption, there exists a constant $\bar \phi > 0$ such that $ \lvert \phi(\mu) \rvert \leq \bar \phi \tilde w(\mu) = \bar \phi \int w d\mu.$ Let $\phi'(\mu) := \phi(\mu) + \bar \phi \int w d\mu$, which is nonnegative. Hence, it is a limit of a nondecreasing sequence of measurable bounded function $\{ \phi'_m \}$ such that $\phi'_m  \uparrow \phi'$. Let $\{ \mu_n \in \mathscr P_w(\mathsf X)\}$ be a converging sequence under $W_1$ to a limit $\mu_0 \in \mathscr P_w(\mathsf X)$. We have then
\begin{align*}
 \liminf_{n \rightarrow \infty} \int \phi'(\mu') \tilde Q(d\mu' | \mu_n, a) \geq \liminf_{n \rightarrow \infty} \int \phi'_m(\mu') \tilde Q(d\mu' | \mu_n, a) = \int \phi'_m(\mu') \tilde Q(d\mu' | \mu_0, a).
\end{align*}
Hence, letting $m \rightarrow \infty$, monotone convergence yields that
\begin{align}
 \liminf_{n \rightarrow \infty} \int \phi'(\mu') \tilde Q(d\mu' | \mu_n, a) \geq \int \phi'(\mu') \tilde Q(d\mu' | \mu_0, a). \label{eq:phi'}
\end{align}
On the other hand, we have for each $n \in \mathbb N$,
\begin{align*}
 \int \tilde w(\mu') \tilde Q(d\mu' | \mu_n, a) = & \int_{\mathsf Y} \int_{\mathsf X} w(x') \tilde P(dx'|\mu_n, a) Q(dy|x', a) = \int_{\mathsf X} \int_{\mathsf X} w(x') P(dx'|x,a) \mu_n(dx).
\end{align*}
Note that by Assumption \ref{asmp:P}(i), we have for each $a\in \mathsf A$, $w'(x, a) := \int_{\mathsf X} w(x') P(dx'|x,a) \in \mathscr L_w$. Proposition \ref{prop:weakconv}(ii) yields that
$
 \lim_{n \rightarrow \infty} \int \tilde w(\mu') \tilde Q(d\mu' | \mu_n, a) = \int \tilde w(\mu') \tilde Q(d\mu' | \mu_0, a).
$
Hence, \eqref{eq:phi'} implies that $\liminf_{n \rightarrow \infty} \int \phi(\mu') \tilde Q(d\mu' | \mu_n, a) \geq \int \phi(\mu') \tilde Q(d\mu' | \mu_0, a).$ In other words, $\mu \mapsto \int \phi(\mu') \tilde Q(d\mu' | \mu, a)$ is lower semicontinuous. We apply this fact to $-\phi$ in lieu of $\phi$ and obtain that $\mu \mapsto \int \phi(\mu') \tilde Q(d\mu' | \mu, a)$ is also upper semicontinuous. Thus the required continuity holds. \proofend

%\begin{lemma}
%Suppose Assumption [REF] holds. If $\phi$ satisfies that $\lvert \phi(\mu)\rvert \leq \frac{\bar r}{1-\alpha \beta} \tilde w(\mu), \forall \mu \in \mathscr P_w$, then $\lvert T(\phi)(\mu)\rvert \leq \frac{\bar r}{1-\alpha \beta} \tilde w(\mu), \forall \mu \in \mathscr P_w$.    
%\end{lemma}

We immediately obtain the following result.

\begin{theorem}
 Suppose Assumptions \ref{asmp:mucontinuity}, \ref{asmp:reward}(i) and \ref{asmp:P}(i) hold. If $\phi \in \mathscr B_{\tilde w}(\mathsf X)$ is convex, then $\mathcal T(\phi)$ is convex and lower semicontinuous.
\end{theorem}

\begin{remark}\label{rm:ascont}
 \citet[Theorem 3.6]{feinberg2016partially} show that one sufficient condition to guarantee both Assumption \ref{asmp:PQ} and \ref{asmp:mucontinuity} is that (i) The stochastic kernel $P(dx'|x,a)$ is canonically weakly continuous and (ii) the stochastic kernel $Q(dy|x,a)$ is continuous in total variation. In addition, it is demonstrated in \cite[Example 4.1]{feinberg2016partially} that the latter continuity in total variation cannot be weakened to the canonical weak continuity. This confirms the necessity of Assumption \ref{asmp:mucontinuity}(ii). 
\end{remark}

%{\bf Remarks on Assumption \ref{asmp:PQ} and \ref{asmp:mucontinuity}.} 

\subsection{Set iteration} Recall that Corollary \ref{coro:fm} imply that a convex and lower semicontinuous function $\phi$ admits a representation of $\phi(\mu) = \sup_{f \in \mathcal N} \int f d\mu$ with some set $\mathcal N \subset \mathscr L(\mathsf X)$. Hence, instead of iterating the value function, we can iterate the acceptance set, which is described as follows. 

\begin{algorithm}\label{algo:1}
 Start with any set $\bar{\mathcal N}_0 \subset \mathscr L(\mathsf X)$. At time $t$, update the acceptance set using the following two steps:
\begin{align*}
 \phi_{t+1}(\mu) = & \max_{a \in \mathsf A} \left( \tilde{r}(\mu, a) + \alpha \int \left( \sup_{f \in \bar{\mathcal N}_t} \int f(x') M(dx'|\mu, a, y) \right) \tilde R(dy|\mu, a) \right)\\
 \bar{\mathcal N}_{t+1} = & \left\{ f\in\mathscr L(\mathsf X) \ \middle| \ \phi_{t+1}(\mu) \geq \int f d\mu , \forall \mu \in \mathscr P_w(\mathsf X) \right\}.
\end{align*}
These steps will be repeated until some stopping criterion is satisfied.
\end{algorithm}
 
An iteration of null level-sets can be analogously designed as above and is therefore omitted. Note that in course of iteration, Lemma \ref{lm:convex} and \ref{lm:lsc} guarantee that $\phi_t$ is convex and lower semicontinuous for each $t = 1, 2, \ldots.$ Hence, Corollary \ref{coro:fm} ensures that $\phi_t(\mu) = \sup_{f \in \bar{\mathcal N}_t} \int f d\mu, \forall \mu \in \mathscr P_w(\mathsf X),$ and for each $t = 1, 2, \ldots.$ By Theorem \ref{thm:vi}(a), we immediately obtain the following result.

\begin{theorem}\label{Secondbigtheorem}
 Suppose Assumptions \ref{asmp:reward}, \ref{asmp:P}, \ref{asmp:PQ}, \ref{asmp:A} and \ref{asmp:mucontinuity} hold. Let $\phi^*$ be the optimal value function for the POMDP, $\bar r > 0$ and $\gamma \in (0,1)$ be the constants as in Theorem \ref{thm:vi}(a). Then, 
 \begin{align*}
  \lVert \phi^* - \phi_t \rVert_{\tilde w} \leq \bar{r} \gamma^t/(1-\gamma), \textrm{ where }  \phi_t(\mu) = \sup_{f \in \bar{\mathcal N}_t} \int f d \mu, \forall \mu \in \mathscr P_w(\mathsf X).
 \end{align*}
\end{theorem}
This implies that the optimal value function $\phi^*$ can be arbitrarily well approximated by a convex and lower semicontinuous function $\phi$ of the dual form. 
\begin{corollary}
 Suppose Assumptions \ref{asmp:reward}, \ref{asmp:P}, \ref{asmp:PQ}, \ref{asmp:A} and \ref{asmp:mucontinuity} hold. For any $\epsilon > 0$, there exists a set $\mathcal N^\epsilon \subset \mathscr L (\mathsf X)$ satisfying
 $
  \lVert \phi^* - \phi^\epsilon \rVert_{\tilde w} \leq \epsilon, \textrm{ where }  \phi^\epsilon(\mu) := \sup_{f \in \mathcal N^\epsilon} \int f d \mu, \forall \mu \in \mathscr P_w(\mathsf X).
 $
\end{corollary}

%Similar to Q-value iteration in MDP literature (see e.g.\ \citet{hernandez1996discrete}), we can iterate acceptance sets depending on $a$.

%\begin{algorithm}\label{algo:2}
% Start with an initial set $\bar{\mathcal N}_0 \subset \mathscr L$ and set $\bar{\mathcal N}_0^a = \bar{\mathcal N}_0, \forall a \in \mathsf A$. At time $t$, update the null level set using the following two steps:
%\begin{align*}
%\phi_{t+1}(\mu, a) = & \tilde r(\mu, a) + \alpha  \int \left( \sup_{f \in \bigcup_{a \in \mathsf A} \bar{\mathcal N}_t^a} \int f(x') M(dx'|\mu, a, y) \right) \tilde R(dy|\mu, a) \\ 
% \bar{\mathcal N}_{t+1}^a = & \left\{ f\in\mathscr L \ \middle| \   \phi_{t+1}(\mu, a) \geq \int f d\mu, \forall \mu \in \mathscr P_w \right\}.
%\end{align*}
%These steps will be repeated until some stopping criterion is satisfied.
%\end{algorithm}

%The above algorithm is also call a \emph{point-based algorithm} (see \citet{pineau2006anytime}) in the literature of POMDPs with applications in robotics, where finite spaces are considered. It is easy to verify the following convergence result. 

%Under Assumptions \ref{asmp:reward}, \ref{asmp:P}, \ref{asmp:PQ}, \ref{asmp:A}, the iteration stated in Algorithm \ref{algo:2} satisfies
%\begin{align*}
%  \lVert \phi^*(\cdot) - \max_{a \in \mathsf A} \phi_t(\cdot, a) \rVert_{\tilde w} \leq \bar{r} \gamma^t/(1-\gamma).
% \end{align*}

\paragraph{A special case: $Q$ is supported by a reference measure} Let us assume that there exists a reference (probability) measure $\varphi$ on $\mathsf Y$ such that $Q(\cdot | x',a) \ll \varphi(\cdot)$ for all $(x',a) \in \mathsf X \times \mathsf A$. Note that POMDPs in many applications satisfy this assumption. For example, the assumption holds automatically if the observation space is finite. Hence, the density of $Q$ w.r.t.\ $\varphi$ exists and is denoted by $q(y|x',a)$. In this case, the iteration can be further simplified. Indeed, it is easy to verify that
$$
M(dx'|\mu, a, y) = \frac{\tilde P(dx'|\mu, a)q(y|x',a)}{\int_{\mathsf X} \tilde P(dx'|\mu, a)q(y|x',a)} \quad \textrm{and} \quad \tilde R(dy|\mu, a) = \int_{\mathsf X} \tilde P(dx'|\mu, a)q(y|x',a) \varphi(dy)
$$
satisfy \eqref{eq:RM}. Under this setup, the calculation of iteration becomes much simpler. Suppose $\phi \in \mathscr B_{\tilde w}$ is convex and lower semicontinuous, then we have by Corollary \ref{coro:fm}, 
\begin{align}
 \mathcal T_{a}(\phi)(\mu) =  \tilde{r}(\mu, a) + \alpha \int \left( \sup_{f \in \bar{\mathcal N}_\phi} \int f(x') \tilde P(dx'|\mu, a)q(y|x',a) \right) \varphi(dy). \label{eq:Ta}
\end{align}
In particular, the continuity of $Q$ in total variation mentioned in Remark \ref{rm:ascont} is
\begin{align*}
  \int \lvert q(y|x_n, a_n) - q(y|x_0, a_0) \rvert \varphi(dy) \rightarrow 0, \textrm{ as } (x_n, a_n) \rightarrow (x_0, a_0).
\end{align*}
\section{Appendix}
\begin{lemma}\label{firstvar}
	For $\mu\in\mathcal{P}_{1}(\mathsf X)$, it holds that
\begin{equation}
	\log\int_{\mathsf X}e^{g}d\mu=\sup_{\nu\in\mathcal{P}_{1}(\mathsf X)}\left\{\int_{\mathsf X}gd\nu - R(\nu|\mu)\right\}, \hspace{16pt}\forall g\in\mathscr L (\mathsf X)
	\end{equation}
\end{lemma}
\begin{proof}
	Let $\mu_{n}$ with $\frac{d\mu_{n}}{d\mu}=e^{\min(g,n)}\cdot\frac{1}{\int_{\mathsf{X}}e^{\min(g,n)}d\mu}.$ For an arbitrary $\nu\in\mathcal{P}_{1}(\mathsf{X})$ with $\mathcal{R}(\nu|\mu)<\infty,$ we have that $\nu$ is absolutely continuous with respect to $\mu_{n}$ and therefore
	\begin{equation}\label{gamw}
	\begin{split}
	\int_{\mathsf X}gd\nu - \mathcal{R}(\nu|\mu)&=\int_{\mathsf X}gd\nu - \int_{\mathsf X}\log\left(\frac{d\nu}{d\mu}\right)d\nu \\
	&=\int_{\mathsf X}gd\nu -\int_{\mathsf X}\log\left(\frac{d\nu}{d\mu_{n}}\right)d\nu - \int_{\mathsf X}\log\left(\frac{d\mu_{n}}{d\mu}\right)d\nu 
	\\&=\int_{\mathsf X}gd\nu-\mathcal{R}(\nu|\mu_{n})+\log\left(\int_{\mathsf{X}}e^{\min(g,n)}d\mu\right)  -\int_{\mathsf X}\min(g,n)d\nu.\end{split}
	\end{equation}
	Now by using the positivity of $\mathcal{R}$ and applying the monotone convergence theorem we get
	\begin{equation*}
\int_{\mathsf X}gd\nu-	\mathcal{R}(\nu|\mu)\leq 	\log\int_{\mathsf X}e^{g}d\mu.
	\end{equation*}
If $\mathcal{R}(\nu|\mu)=\infty,$ the above inequality holds trivially.\newline

Now, on the other hand, by setting $\nu=\mu_{n}$ in \eqref{gamw}, observing that $\int_{\mathsf X}gd\mu_{n}-\int_{\mathsf X}\min(g,n)d\mu_{n}\geq 0,$ and applying the monotone convergence theorem one more time, we get

\begin{equation*}
\lim_{n\rightarrow \infty}\left(\int_{\mathsf X}gd\mu_{n} - \mathcal{R}(\mu_{n}|\mu) \right)\geq \log\int_{\mathsf X}e^{g}d\mu,
\end{equation*}
which yields our result.
\end{proof}
\begin{lemma}\label{con=lip}
	For $\mu,\nu\in\mathcal{P}_{1}(\mathsf X)$, it holds that
	\begin{equation}
\sup_{g\in C_{b}(\mathsf X)}\left\{\int_{\mathsf X}gd\nu-\log\int_{\mathsf X}e^{g}d\mu\right\}=\sup_{g\in \mathscr L(\mathsf X)}\left\{\int_{\mathsf X}gd\nu-\log\int_{\mathsf X}e^{g}d\mu\right\}
	\end{equation}
\begin{proof}
	For simplicity we will set $F(g)=\int_{\mathsf X}gd\nu-\log\int_{\mathsf X}e^{g}d\mu.$ By properties of the supermum, we have
	\begin{equation*}
	\sup_{g\in C_{b}(\mathsf X)\cap \mathscr L(\mathsf X)}\hspace{-12pt}F(g)\leq 	\sup_{g\in C_{b}(\mathsf X)}\hspace{-5pt}F(g)\leq\hspace{-7pt} 	\sup_{g\in C_{b}(\mathsf X)\cup \mathscr L(\mathsf X)}\hspace{-12pt}F(g), \hspace{8pt}	\sup_{g\in C_{b}(\mathsf X)\cap \mathscr L(\mathsf X)}\hspace{-12pt}F(g)\leq	\sup_{g\in  \mathscr L(\mathsf X)}\hspace{-5pt}F(g)\leq \hspace{-7pt}	\sup_{g\in C_{b}(\mathsf X)\cup \mathscr L(\mathsf X)}\hspace{-12pt}F(g). 
	\end{equation*}
	
    So it will be enough to prove that
    \begin{equation*}
    	\sup_{g\in C_{b}(\mathsf X)\cup \mathscr L(\mathsf X)}\hspace{-12pt}F(g)\hspace{4pt}\leq \sup_{g\in C_{b}(\mathsf X)\cap \mathscr L(\mathsf X)}\hspace{-12pt}F(g).
    \end{equation*}
    Let $g\in C_{b}(\mathsf X)\cup \mathscr L(\mathsf X),$ with $F(g)\neq -\infty.$  It is now enough to prove that for every $ \epsilon>0,$ it exists $\tilde{g}\in C_{b}(\mathsf X)\cap \mathscr L(\mathsf X)$ such that $F(g)-F(\tilde{g})\leq \epsilon.$  
    
    First we are going to approximate $g$ by a bounded function $\hat{g}.$ We set $g_{m}= min(g,m),$ and by a suitable application of the monotone convergence theorem we have that for sufficiently big $m$ holds $|F(g)-F(g_{m})|\leq\frac{\epsilon}{4}.$ Now if we further set $g^{n}_{m}=\max(g_{m},-n),$ and we apply the dominated convergence theorem, we can find $\hat{g}=g^{n}_{m}$ such that 
    
    \begin{equation}\label{key}
    |F(g)-F(\hat{g})|\leq \frac{\epsilon}{2}.
    \end{equation}
    
    To get the Lipschiz property, we will first apply Prohorov's theorem, and we will find a compact set $\mathsf K$ such that $\mu(\mathsf X\setminus \mathsf K) ,\nu(\mathsf X\setminus \mathsf K) \leq \epsilon'.$  Now, we can approximate any function $\hat{g}$ in  $C_{b}(\mathsf K),$ $\epsilon'$-uniformly by a function in $\mathscr L(\mathsf K)$ through the formula 
	\begin{equation}
	\hat{g}_{n}(x)=\inf_{y\in\mathsf{K}}\left\{\hat{g}(y)+ nd(x,y)\right\}, 
	\end{equation}
	for sufficiently large $n$. By taking $n$ big enough we can also have that $e^{\hat{g}_{n}},e^{\hat{g}},$ are at most $\epsilon'$-uniformly apart in $\mathsf{K}.$  This formula actually defines $\hat{g}_{n}$ to be Lipschitz on the whole space $\mathsf{X}.$ We can further bound by using $\tilde{g}(x)=\max\left(\min\left(\hat{g}_{n}(x),\|\hat{g}\|_{\infty}\right),-\|\hat{g}\|_{\infty}\right).$
Now we have 
\begin{equation}
\begin{split}
\left|\int_{\mathsf{K}}\hat{g}d\nu-\int_{\mathsf{K}}\tilde{g}d\nu\right|\leq \epsilon',\hspace{8pt} \left|\int_{\mathsf X\setminus\mathsf{K}}\hat{g}d\nu-\int_{\mathsf X\setminus\mathsf{K}}\tilde{g}d\nu\right|\leq 2\epsilon'\|\hat{g}\|_{\infty}\\ \left|\int_{\mathsf{K}}e^{\hat{g}}d\nu-\int_{\mathsf{K}}e^{\tilde{g}}d\nu\right|\leq \epsilon',\hspace{8pt} \left|\int_{\mathsf X\setminus\mathsf{K}}e^{\hat{g}}d\nu-\int_{\mathsf X\setminus\mathsf{K}}e^{\tilde{g}}d\nu\right|\leq 2\epsilon' e^{\|\hat{g}\|_{\infty}}
\end{split}
\end{equation}
Now by using the modulus of uniform continuity $\omega$ for the logarithm on $[e^{-\|\hat{g}\|_{\infty}},e^{\|\hat{g}\|_{\infty}}].$ (or a simple mean value theorem), we get the following estimate 
$|F(\hat{g})-F(\tilde{g})|\leq 2\epsilon' +2\epsilon'\|\hat{g}\|_{\infty}+\omega(2\epsilon' (1+e^{\|\hat{g}\|_{\infty}})). 
$
Now if $\epsilon'$ becomes sufficiently small we have $|F(\hat{g})-F(\tilde{g})|\leq\frac{\epsilon}{2},$ and by combining  with
\eqref{key} we get our claim.
\end{proof}
\end{lemma}

We conclude by providing a proof for Corollary \ref{th:vil}.

\begin{proof}[Proof of Corollary \ref{th:vil}]
	First assume that $\Phi( W_1(\mu, \nu) )\leq \phi(\mu).$ For $f\in [L(\mathsf X)]_{1},$ we have
	\begin{equation*}
	\begin{split}
	\rho\left(t\int_{X}fd\nu-tf-\Phi^{*}(t)\right)=&\sup_{\mu\in\mathcal{P}_{1}(\mathsf{X})}\left[\int_{\mathsf{X}}\left(t\int_{\mathsf{X}}f d\nu-tf-\Phi^{*}(t)\right)d\mu-\phi(\mu)\right]\\=&\sup_{\mu\in\mathcal{P}_{1}(\mathsf{X})}\left[t\left(\int_{\mathsf{X}}
	fd\nu-\int_{\mathsf{X}}fd\mu\right)-\Phi^{*}(t)-\phi(\mu)\right]\\=&
	\sup_{\mu\in\mathcal{P}_{1}(\mathsf{X})}\left[tW_{1}(\mu,\nu)-\Phi^{*}(t)-\phi(\mu)\right]
	\\=&\sup_{\mu\in\mathcal{P}_{1}(\mathsf{X})}\left[\Phi(W_{1}(\mu,\nu))-\phi(\mu)\right]\leq 0.		\end{split}
	\end{equation*}
	Conversely if 	$\rho\left(t\int_{X}fd\nu-tf-\Phi^{*}(t)\right)\leq 0,$
	we have
	\begin{equation}
	\begin{split}
	t\left(\int_{\mathsf{X}}
	fd\nu-\int_{\mathsf{X}}fd\mu\right)-\Phi^{*}(t)=&\int_{\mathsf{X}}\left(t\int_{\mathsf{X}}f d\nu-tf-\Phi^{*}(t)\right)d\mu\\
	&\leq \rho\left(t\int_{X}fd\nu-tf-\Phi^{*}(t)\right) +\phi(\mu)\leq\phi(\mu).	
\end{split}
	\end{equation}
	Taking the supremum over $f\in [\mathscr L(\mathsf X)]_{1},$ we get 
	\begin{equation}
	tW_{1}(\mu,\nu)-\Phi^{*}(t)\leq \phi(\mu)
	\end{equation}
By taking the supremum over $t\geq 0,$ we have  $\Phi( W_1(\mu, \nu) )\leq \phi(\mu).$
\end{proof}

\end{document}